\documentclass[12pt]{amsart}

\usepackage{graphicx}
\usepackage{amssymb}
\usepackage{amsmath}
\usepackage{mathtools}
\usepackage{seqsplit}
\usepackage[utf8]{inputenc}
\usepackage{float}
\usepackage{longtable}
\usepackage{url}
\usepackage{enumitem}
\usepackage{indentfirst}
\usepackage{bbm}
\usepackage{xcolor}

\usepackage{lineno}
\usepackage{amsthm}
\usepackage{dsfont}
\usepackage{hyperref}
\usepackage{tikz}
\usepackage[left=2.5cm, right=2.5cm, top=2.5cm, bottom=2.5cm]{geometry}
\usetikzlibrary{positioning}
  \def\gn#1#2{{$\href{http://groupnames.org/\#?#1}{#2}$}}
\def\gn#1#2{$#2$}  
\tikzset{sgplattice/.style={inner sep=1pt,norm/.style={red!50!blue},char/.style={blue!50!black},
  lin/.style={black!50}},cnj/.style={black!50,yshift=-2.5pt,left=-1pt of #1,scale=0.5,fill=white}}
\usetikzlibrary{decorations.pathmorphing,calc,shapes.misc}
\usepackage{cleveref}

\newcommand{\Diff}{\operatorname{Diff}}

\newcommand{\pt}{\operatorname{pt}}

\newcommand{\SL}{\operatorname{SL}}
\newcommand{\TL}{\operatorname{TL}}

\newcommand{\chiV}{\chi_{\raisebox{-0.2ex}{$\scriptstyle V$}}}
\newcommand{\chiCk}{\chi_{\raisebox{-0.2ex}{$\scriptstyle\mathbb{C}_k$}}}

\newtheorem{theorem}{Theorem}[section]
\newtheorem{definition}[theorem]{Definition}
\newtheorem{lemma}[theorem]{Lemma}
\newtheorem{corollary}[theorem]{Corollary}

\newtheorem{remark}[theorem]{Remark}

\newtheorem*{remark*}{Remark}

\newtheorem*{theorem*}{Theorem}
\newtheorem*{conjecture*}{Conjecture}
\newtheorem*{question*}{Question}

\numberwithin{table}{section}
\numberwithin{figure}{section}

\makeatletter
\def\ps@pprintTitle{%
  \let\@oddhead\@empty
  \let\@evenhead\@empty
  \def\@oddfoot{}%
  \let\@evenfoot\@oddfoot
}
\makeatother

\date{}  

\title[Nonexistence of one-fixed-point actions of direct products on spheres]{Nonexistence of one-fixed-point actions on spheres of direct products possessing a cyclic group as a factor}

\author{Piotr Mizerka}
\author{Kamil Śmietaniak}


\makeatletter
\AtBeginDocument{%

\renewenvironment{proof}[1][\proofname]{%
  \par\pushQED{\qed}%
  \normalfont
  \topsep6\p@\@plus6\p@\relax
  \trivlist
  \item[\hskip\labelsep\normalfont #1\@addpunct{.}]\ignorespaces
}{%
  \popQED\endtrivlist\@endpefalse
}

}%
\makeatother

\begin{document}

\begin{abstract}
In this paper, we show that there is no smooth, effective, one–fixed-point action of 
\( G \in \{ \mathrm{SL}(2,5) \times C_k,\, \mathrm{TL}(2,5) \times C_k,\, S_5 \times C_k \} \) 
on \( S^n \) for certain dimensions \( n \ge 6 \).
For 
\( G \in \{ \mathrm{SL}(2,5) \times C_k,\, \mathrm{TL}(2,5) \times C_k \} \) 
with \( k \) coprime to \( 30 \), 
this result completes the classification of dimensions \( n \) 
for which \( G \) admits a smooth, effective, one–fixed-point action on \( S^n \).

\end{abstract}

\maketitle







\setcounter{section}{-1}

\section{Introduction}
This article is concerned with smooth actions of finite groups on manifolds. Unless specified otherwise, all groups are assumed to be finite. A \emph{smooth effective action} of a group $G$ on a manifold $M$ is a monomorphism $G \rightarrowtail \Diff(M)$, where $\Diff(M)$ denotes the group of all diffeomorphisms of $M$. In our setting, $M$ is a sphere and $G$ is a direct product of a specific finite group and a cyclic group. Unless stated otherwise, we will henceforth assume that all actions under consideration are smooth and effective.

It is well known which groups admit a one–fixed-point action on a sphere \cite{Oliver1974,Oliver1975,Oliver1976,Oliver1977,Oliver1978,Morimoto1998}. Such groups $G$ are precisely the \emph{Oliver groups}, characterized by the following property \cite{Oliver1974,Oliver1975,Oliver1977,Oliver1978}: there does not exist a sequence of subgroups $P \trianglelefteq H \trianglelefteq G$ such that $P$ and $G/H$ are groups of prime power order and $H/P$ is cyclic. It turns out that this condition is even necessary for the existence of a one fixed point action on a homology sphere. Indeed, suppose that there exists a sequence of normal subgroups $P\trianglelefteq H\trianglelefteq G$, where $P$ is a $p$-group, $H/P$ is cyclic, and $G/H$ is a $q$-group for some primes $p$ and $q$. Notice that if a one fixed point action on a homology sphere $\Sigma$ existed, then by the so called Slice Lemma (\Cref{theorem:slice}), we would obtain a fixed point free action of $G$ on a $\mathbb{Z}$-acyclic complex. We claim that $G$ cannot act without fixed points on a finite $\mathbb{Z}$-acyclic $G$-CW complex $X$, and hence, in particular, on a compact $\mathbb{Z}$-acyclic manifold. Indeed, the argument proceeds as follows.
 
\begin{enumerate}
\item By Smith theory, the fixed point set $X^{P}$ is
$\mathbb{Z}_{p}$-acyclic. In particular, $X^{P}$ is non-empty and
rationally acyclic.
 
\item Since $P\trianglelefteq H$, the cyclic group $H/P$ acts on
$X^{P}$. Let $\overline{h}$ be a generator of $H/P$. Then
\[
    X^{H}=(X^{P})^{H/P}=(X^{P})^{\overline{h}}.
\]
By the Lefschetz fixed point formula for finite-order cellular
transformations,
\[
    \chi(X^{H})
    =
    \chi\bigl((X^{P})^{\overline{h}}\bigr)
    =
    L(\overline{h}),
\]
where $\overline{h}\colon X^P\rightarrow X^P$. Since \(X^{P}\) is rationally acyclic, the only non-zero rational
homology group of \(X^{P}\) is \(H_{0}(X^{P};\mathbb{Q})\cong\mathbb{Q}\),
on which \(\overline{h}\) acts trivially. Consequently,
\[
    L(\overline{h})=1,
\]
and therefore
\[
    \chi(X^{H})=1.
\]
 
\item Finally, since \(H\trianglelefteq G\), the \(q\)-group \(G/H\)
acts on \(X^{H}\). Counting the cells of \(X^{H}\) in \(G/H\)-orbits
shows that
\[
    \chi(X^{H})
    \equiv
    \chi\bigl((X^{H})^{G/H}\bigr)
    =
    \chi(X^{G})
    \pmod q.
\]
If the action of \(G\) were fixed-point-free, then \(X^{G}\) would be
empty and hence \(\chi(X^{G})=0\). This would give
\[
    1=\chi(X^{H})\equiv\chi(X^{G})=0\pmod q,
\]
which is impossible.
 
\end{enumerate} 
Thus, the existence of such a normal series forces every action of $G$ on a
finite $\mathbb{Z}$-acyclic complex to have a global fixed point.

Recently, attention has been focused on determining, for a given Oliver group $G$, the set of all dimensions $\mathbb{N}_G$ of spheres on which $G$ admits an action with exactly one fixed point. It is known that $\mathbb{N}_G \subseteq \{6,7,\ldots\}$ \cite{Morimoto1987,Furuta1989,DeMichelis1989,Schultz1990}, and that $6 \in \mathbb{N}_G$ implies $G \in \{A_5, A_5 \times C_2, S_5\}$ \cite{Tamura2025}, where $A_5$ and $S_5$ denote the alternating and symmetric groups of degree $5$, respectively. The first group $G$ for which the list $\mathbb{N}_G$ was determined is $A_5$ \cite{Morimoto1987,BakMorimoto,mornonslvb1}. Historically, the first construction of a one–fixed-point action of a finite group on a sphere was that of $\SL(2,5)$ — the group of $2 \times 2$ matrices over $\mathbb{F}_5$ with determinant one — acting on $S^7$ \cite{stein}. In contrast, it is known that $8 \notin \mathbb{N}_{\SL(2,5)}$ \cite{Borowiecka2016}. Further exclusions were obtained in \cite{boromiz}, as well as in \cite{tamuraventi} and \cite{ogor2023,morimoto_appendix_mizerka}. Finally, in \cite{Morimoto2025_OneFixedPointIII} (preceded by \cite{mornonslvb1} and \cite{morimotononsol2}), the sets $\mathbb{N}_G$ were computed for $\SL(2,5)$, $\TL(2,5)$, and $A_5 \times C_p$ with $p$ prime to $30$, where $\TL(2,5)$ denotes the unique non-split extension of $C_2$ by $\SL(2,5)$. Moreover, \cite{Morimoto2025_OneFixedPointIII} provides constructions of one–fixed–point actions not only for the cases in which the sets $\mathbb{N}_G$ were completely determined, but also for groups of the form $G = H \times C_p$, where $H \in \{\SL(2,5), \TL(2,5)\}$.

In this paper, we establish the nonexistence of one–fixed–point actions of 
$G \in \{\SL(2,5) \times C_k,\, \TL(2,5) \times C_k,\, S_5 \times C_k\}$ 
on spheres in certain dimensions. For the cases 
$G \in \{\SL(2,5) \times C_k,\, \TL(2,5) \times C_k\}$ 
with $k$ coprime to $30$, these results, together with the current state of knowledge, allow us to determine the corresponding sets $\mathbb{N}_G$. The main results of the paper are as follows.

\begin{theorem}\label{theorem:1main}\emph{(cf. \Cref{theorem:1main2}, \Cref{theorem:1main2only}, \Cref{theorem:1main2_even})}
    Let $G = \SL(2,5) \times C_k$ with $k$ odd. Then there does not exist a one–fixed–point action of $G$ on $S^n$ for any $n \in \{8,9\}$. In the case where $k$ is even, there does not exist a one–fixed–point action of $G$ on $S^n$ for any $n \in \{7,8,9\}$.
\end{theorem}

\noindent As a consequence of the theorem above, together with the results of \cite{Tamura2025} and \cite[Theorem~4.1]{Morimoto2025_OneFixedPointIII}, we obtain
\[
    \mathbb{N}_{\SL(2,5) \times C_k} = \{7\} \cup \{10,11,\ldots\}
\]
for $k$ coprime to $30$.

\begin{theorem}\label{theorem:2main}\emph{(cf. \Cref{theorem:2main2}, \Cref{theorem:2main2only}, \Cref{theorem:2main2_even})}
    Let $G = \TL(2,5) \times C_k$ with $k$ odd. Then there does not exist a one–fixed–point action of $G$ on $S^n$ for any $n \in \{6,\ldots,13\} \cup \{15,16,17,21\}$. In the case where $k$ is even, there does not exist a one–fixed–point action of $G$ on $S^n$ for any $n \in \{6,\ldots,17\} \cup \{21\}$.
\end{theorem}

\noindent As a consequence of the theorem above and \cite[Theorem~4.3]{Morimoto2025_OneFixedPointIII}, we obtain
\[
    \mathbb{N}_{\TL(2,5) \times C_k} = \{14,18,19,20,22,23,\ldots\}
\]
for $k$ coprime to $30$.

\begin{theorem}\label{theorem:3main}\emph{(cf. \Cref{theorem:3main2}, \Cref{lemma:3_8}, \Cref{theorem:3main2only}, \Cref{theorem:3main22})}
    Let $G = S_5 \times C_k$ with $k > 2$. Then there does not exist a one–fixed–point action of $G$ on $S^n$ for any $n \in \{6,\ldots,13\} \cup \{15\}$ in the case $k$ is odd and for any $n \in \{6,\ldots,13\} \cup \{15,17\}$ in the case $k$ is even or $k$ is divisible by $3$. In the case where $k = 2$, there does not exist a one–fixed–point action of $G$ on $S^n$ for any $n \in \{6,\ldots,9\} \cup \{13\}$.
\end{theorem}

This paper is organized as follows. In \cref{section:2}, we provide the fixed point dimensions for the irreducible representations of the groups under consideration, together with their posets of conjugacy classes of subgroups. In \cref{section:1}, we present all the necessary exclusion techniques used to obtain our results. These consist primarily of previously established methods, supplemented by several new ones. \Cref{section:3} forms the central part of this article, containing the proofs of \Cref{theorem:1main}, \Cref{theorem:2main}, and \Cref{theorem:3main}. In the appendices, we also present exclusion results for $A_5 \times C_k$, as well as for extensions of finite groups by a torus.

\section*{Acknowledgements}
We would like to thank Prof. Masaharu Morimoto for his valuable suggestions and insights regarding actions of $A_5$. We are also grateful to Prof. Krzysztof Pawałowski for his support and many helpful remarks during the numerous seminar meetings dedicated to this work. We thank Dr Shunsuke Tamura for his interest and useful remarks, which substantially improved the presentation of this paper.

We additionally acknowledge the use of \textit{ChatGPT} \cite{ChatGPT2025} as a writing assistant for polishing the English language and improving the \LaTeX{} presentation of the manuscript. The whole contents was also verified with \textit{ChatGPT} \cite{ChatGPT2025} and a few errors were spotted and corrected accordingly.

\section{Subgroup structures and fixed point dimension tables}\label{section:2}
We recall first how to get real irreducible representations from the complex ones. The contents of the discussion below is copied verbatim from the PhD thesis of the first author \cite[subsection 2.1.2]{Mizerka2021}.

Let $V$ be an irreducible $\mathbb{C}G$-module with character $\chi$. There are three mutually exclusive possibilities (see \cite[Proposition 38]{Serre1977}).
\begin{enumerate}
    \item There is no nonzero invariant bilinear form on $V$. In this case $\chi(g)$ is not real for some $g\in G$. Moreover, the character $2\operatorname{Re}(\chi)=\chi+\overline{\chi}$ is the character of an irreducible $\mathbb{R}G$-module.
    \item There is a nonzero symmetric invariant bilinear form on $V$. Then $\chi$ is \emph{realizable} over $\mathbb{R}$, that is, the $\mathbb{RG}$-module $V$ is isomorphic (as $\mathbb{C}G$-module) to some $\mathbb{R}G$-module.
    \item There exists a nonzero skew-symmetric invariant bilinear form on $V$ -- in such a case, $\chi$ is real but not realizable over $\mathbb{R}$ and $2\chi$ is the character of some irreducible $\mathbb{R}G$-module.
\end{enumerate}
Furthermore, every irreducible $\mathbb{R}G$-module can be obtained in one of the three ways described above, see \cite[pp. 108]{Serre1977}. 

It turns out that for the determination of irreducible $\mathbb{R}G$-modules, the concept of \emph{Frobenius-Schur indicator} is very useful. It is defined below.
\begin{definition}\emph{\cite{Serre1977}}\label{notation:indicator}
The \textbf{Frobenius-Schur indicator} of a $\mathbb{C}G$-module $V$ with character $\chi$ is defined to be
$$
\iota(\chi)=\frac{1}{|G|}\sum_{g\in G}\chi(g^2).
$$
\end{definition}
\noindent For irreducible $\mathbb{C}G$-modules, the Frobenius-Schur indicator can take one of the three values: $0$, $1$ or $-1$, see \cite[pp. 109, Proposition 39]{Serre1977}.
The following theorem allows us to deduce existence of appropriate bilinear forms on irreducible $\mathbb{C}G$-modules with the help of the Frobenius-Schur indicator.
\begin{theorem}\emph{\cite[23.16 Theorem]{Liebeck2001}}
Let $V$ be an irreducible $\mathbb{C}G$-module with character $\chi$. Then, the following statements hold.
\begin{enumerate}
    \item $\iota(\chi)=0$ if and only if there is no nonzero invariant bilinear form on $V$,
    \item $\iota(\chi)=1$ if and only if there exists a nonzero symmetric bilinear form on $V$,
    \item $\iota(\chi)=-1$ if and only if there is a nonzero skew-symmetric bilinear form on $V$.
\end{enumerate}
\end{theorem}
Thus, we may determine characters of irreducible $\mathbb{R}G$-modules from complex ones as follows. We look at all irreducible complex characters $\chi$ and compute their Frobenius-Schur indicators $\iota(\chi)$. These indicators contain information about which one of the characters from $\chi,2\operatorname{Re}(\chi)$ is a character of an irreducible $\mathbb{R}G$-module. The set of so obtained characters (we omit repetitions -- the characters of the form $2\operatorname{Re}(\chi)$ appear twice -- once for the irreducible complex character $\chi$ and second for $\overline{\chi}$ which is irreducible as well) is the complete list of characters of real irreducible representations.
\newline
\newline
\indent We are interested in the groups of the form $G\times C_k$. Their irreducible representations can be grouped into three classes:

\begin{itemize}
    \item \textbf{Type I (pure $G$)} – a $G$-representation with trivial $C_k$-action.
    \item \textbf{Type II (pure $C_k$)} – a $C_k$-representation with trivial $G$-action. These representations have two different dimensions: the one-dimensional $\mathbb{R}_{\pm}$ associated to the subgroup of order two, and the two-dimensional $\mathbb{C}_k$, which arise from the irreducible (not necessarily faithful) complex representations of $C_k$.
    \item \textbf{Type III (mixed)} – the tensor product of a Type I and a Type II representation.
\end{itemize}
\noindent More details on the last type from the trichotomy above are discussed in the paragraph after \Cref{lemma:fixed_point_dims2}.

To compute the dimensions of fixed point subspaces for $H \leq G$, we will use the following well-known lemma.

\begin{lemma}\label{lemma:fixed_point_dims2}
    Let $V$ be a complex $G$-representation with character $\chiV$. Then, for any subgroup $H \leq G$, one has
    \[
        \dim V^H = \frac{1}{|H|} \sum_{h \in H} \chiV(h).
    \]
\end{lemma}

In the case an irreducible real $G$-representation $V$ arises from the complex irreducible $G$-representation with Frobenius-Schur indicator equal to $1$, we write $V \otimes \mathbb{C}_k$ for any of the $G \times C_k$-representations whose restrictions to $G$ and $C_k$ are isomorphic to two copies of $V$ and a $C_k$-module with $a^2$ acting non-trivially, where $a$ is the generator of $C_k$. In the analogous cases when Frobenius-Schur indicator is equal to $-1$ or $0$, we write $V \otimes_{\mathbb{C}} \mathbb{C}_k$. In the case $a^2$ acts trivially, we write $V \otimes \mathbb{R}_{\pm}$ for the $G \times C_k$-representation that acts faithfully on the second component via the quotient $C_2$ (this applies only when $k$ is even). The real irreducible representations of $G\times C_k$ are precisely those of the form $V$, $V\otimes\mathbb{C}_k$, $V\otimes_{\mathbb{C}}\mathbb{C}_k$, and $V\otimes \mathbb{R}_{\pm}$. There is an ambiguity in the introduced definitions, however, which is of no harm to us -- see \Cref{remark:ambig}.

Since we are only interested in the dimensions of fixed point subspaces for the whole group $C_k$, the following lemma provides a complete answer.

\begin{corollary}\label{corollary:fixed_dims2}
Let $H \leq G$, and let $W\in\{V \otimes \mathbb{C}_{k},V\otimes_{\mathbb{C}}\mathbb{C}_k\}$ be an irreducible $G \times C_k$–representation. Then:
\begin{enumerate}
    \item $\dim W^{H \times C_k} = 0$,
    \item $\dim W^{H} = 2 \dim V^H$ if $W=V\otimes\mathbb{C}_k$ and $\dim W^{H} = \dim V^H$ if $W=V\otimes_{\mathbb{C}}\mathbb{C}_k$.
\end{enumerate}
\end{corollary}

\begin{proof}
    For the first case, we compute the dimension directly from the characters of $V$ and $\mathbb{C}_k$, using the fact that the character of a tensor product of representations is the product of the individual characters (the real fixed point set dimension of the real module $W$ is equal to the complex fixed point set dimension of the complexification of $W$ for which we may use the formula from \Cref{lemma:fixed_point_dims2}):
    \[
    |H \times C_k| \dim W^{H\times C_k}
    = \sum_{(h, l) \in H \times C_k} \chiV(h) \chiCk(l)
    = \sum_{h \in H} \chiV(h) \sum_{l \in C_k} \chiCk(l)
    = \sum_{h \in H} \chiV(h) \cdot 0 = 0.
    \]
    The second conclusion follows from the fact that, after forgetting the $C_k$–action, the resulting representation is isomorphic to two copies of $V$ in the real case and to one copy of $V$ in the complex case.
\end{proof}

\begin{remark}\label{remark:ambig}
\emph{Note that, so far, we have slightly abused the notation $\mathbb{C}_k$. The symbol $\mathbb{C}_k$ may denote any representation with trivial $G$–action and non-trivial $C_k$–action for some $a \in C_k$ of order at least $3$. In light of \Cref{corollary:fixed_dims2}, this abuse is justified, as the choice of such a representation does not affect the fixed point dimensions for the subgroups under consideration. The reason we make this simplification is that, in all arguments referring to representations, the fixed point dimensions provide sufficient information.}
\end{remark}

Below, we present the subgroup structures (by conjugacy classes) of $\SL(2,5)$, $\TL(2,5)$, and $S_5$, together with their fixed point dimension tables for non-trivial real irreducible representations. In fact, we are interested only in the fixed point dimensions corresponding to each given representation. Hence, if two representations have the same fixed point dimensions for all subgroups, we identify them. The conjugacy classes of subgroups below and \Cref{fig:A5} were taken from \cite{GroupNames}. To obtain the fixed point dimension tables, including the determination of characters of real irreducible representations, we used our GAP \cite{GAP4} package \texttt{OneFPActions} \cite{mizerka2020one_fixed_point_spheres_scripts}. Details on reproducing these tables are provided in \cite{mizerka2020one_fixed_point_spheres_scripts}. In the posets of subgroup conjugacy classes, each conjugacy class is labeled with the number of subgroups it contains, indicated in the lower left corner. If no such indicator is present, the corresponding subgroup is normal.

\begin{figure}[H]
  \centering
    \begin{tikzpicture}[scale=1.0,sgplattice]
  \node at (2,0) (1) {\gn{C1}{C_1}};
  \node at (2,0.803) (2) {\gn{C2}{C_2}};
  \node at (3.88,1.76) (3) {\gn{C3}{C_3}};
  \node at (0.125,1.76) (4) {\gn{C5}{C_5}};
  \node at (2,1.76) (5) {\gn{C4}{C_4}};
  \node at (3.88,2.97) (6) {\gn{C6}{C_6}};
  \node at (0.125,2.97) (7) {\gn{C10}{C_{10}}};
  \node at (2,2.97) (8) {\gn{Q8}{Q_8}};
  \node at (3.88,4.19) (9) {\gn{Dic3}{{C_{3}\rtimes C_{4}}}};
  \node at (0.125,4.19) (10) {\gn{Dic5}{{C_{5}\rtimes C_{4}}}};
  \node at (2,4.19) (11) {\gn{SL(2,3)}{{\rm SL(2,3)}}};
  \node at (2,5.14) (12) {\gn{SL(2,5)}{{\rm SL(2,5)}}};
  \draw[lin] (1)--(2) (1)--(3) (1)--(4) (2)--(5) (2)--(6) (3)--(6) (2)--(7)
     (4)--(7) (5)--(8) (5)--(9) (6)--(9) (5)--(10) (7)--(10) (6)--(11)
     (8)--(11) (9)--(12) (10)--(12) (11)--(12);
  \node[cnj=3] {10};
  \node[cnj=4] {6};
  \node[cnj=5] {15};
  \node[cnj=6] {10};
  \node[cnj=7] {6};
  \node[cnj=8] {5};
  \node[cnj=9] {10};
  \node[cnj=10] {6};
  \node[cnj=11] {5};
  
\end{tikzpicture}
  \caption{Conjugacy classes of subgroups of $\SL(2,5)$.}
  \label{fig:sl25}
\end{figure}
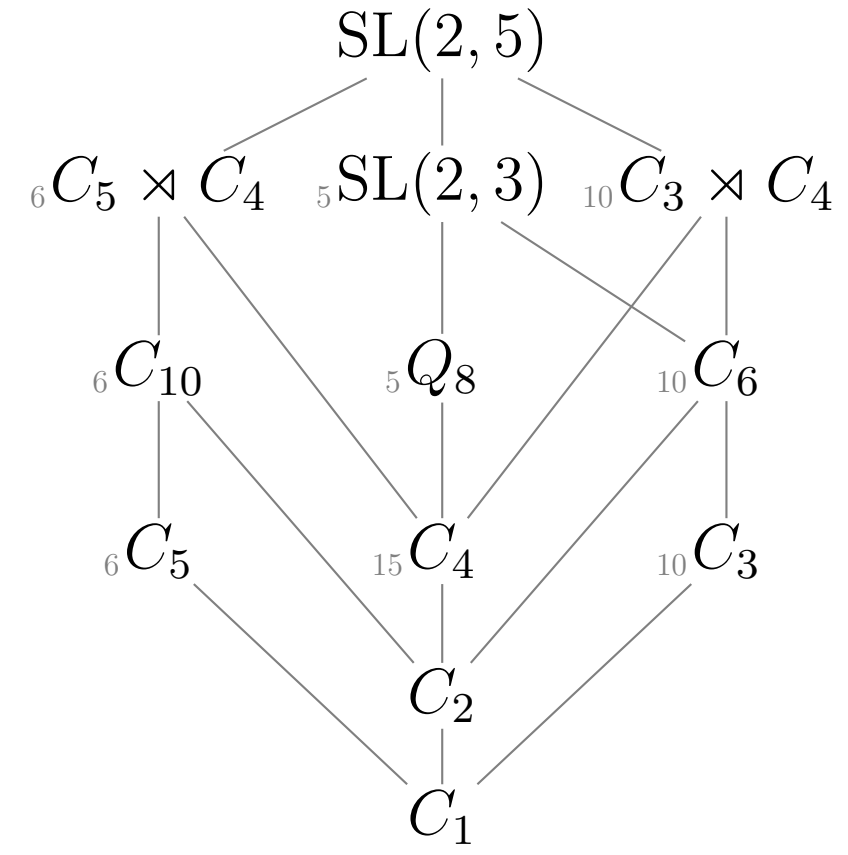

\begin{table}[H]
\centering
{\fontsize{10pt}{12pt}\selectfont 
\begin{tabular}{|c||c|c|c|c|c|c|c|c|c|c|c|c|c|c|c|c|c|c|c|c|c|} \hline
& 
$C_1$ & $C_{2}$  &$C_{3}$ & 
$C_{4}$  &$C_{5}$ & $C_{6}$ & $Q_{8}$ &$C_{10}$& 

 $C_{3}\rtimes C_{4}$ & $C_{5}\rtimes C_{4}$ & $\mathrm{SL}(2,3)$ & $\mathrm{SL}(2,5)$
 \\ 
\hline
\hline

$V_3$ & 
$3$ & $3$ & $1$ & $1$ & $1$ & 
$1$ & $0$ &$1$ & 
 $0$ &$0$ & $0$ & $0$ \\
\hline

$V_{4,1}$ & 
$4$ & $4$ & $2$ & $2$ & $0$ & 
$2$ & $1$ & $0$ & 
 $1$ &$0$ & $1$ & $0$ \\
\hline 

$V_{4,2}$ & 
$4$ & $0$ & $0$ & $0$ & $0$ & 
$0$ & $0$  & $0$ & 
 $0$ &$0$ & $0$ & $0$ \\
\hline

$V_{5}$ & 
$5$ &$5$& $1$ & $3$ & $1$ & $1$ & 
$2$ & $1$ & 
 $1$ &$1$ & $0$ & $0$  \\
\hline 

$V_{8}$ & 
$8$ & $0$ & $4$ & $0$ & $0$ & 
$0$ & $0$  & $0$ & 
 $0$ &$0$ & $0$ & $0$  \\

\hline 

$V_{12}$ & 
$12$ & $0$ & $4$ & $0$ & $4$ & 
$0$ & $0$ & $0$ & 
 $0$ &$0$ & $0$ & $0$ \\
\hline

\end{tabular}
}

\caption{Fixed point dimension table for the non-trivial real irreducible representations of $\SL(2,5)$.}\label{table:dims_sl25}
\end{table}

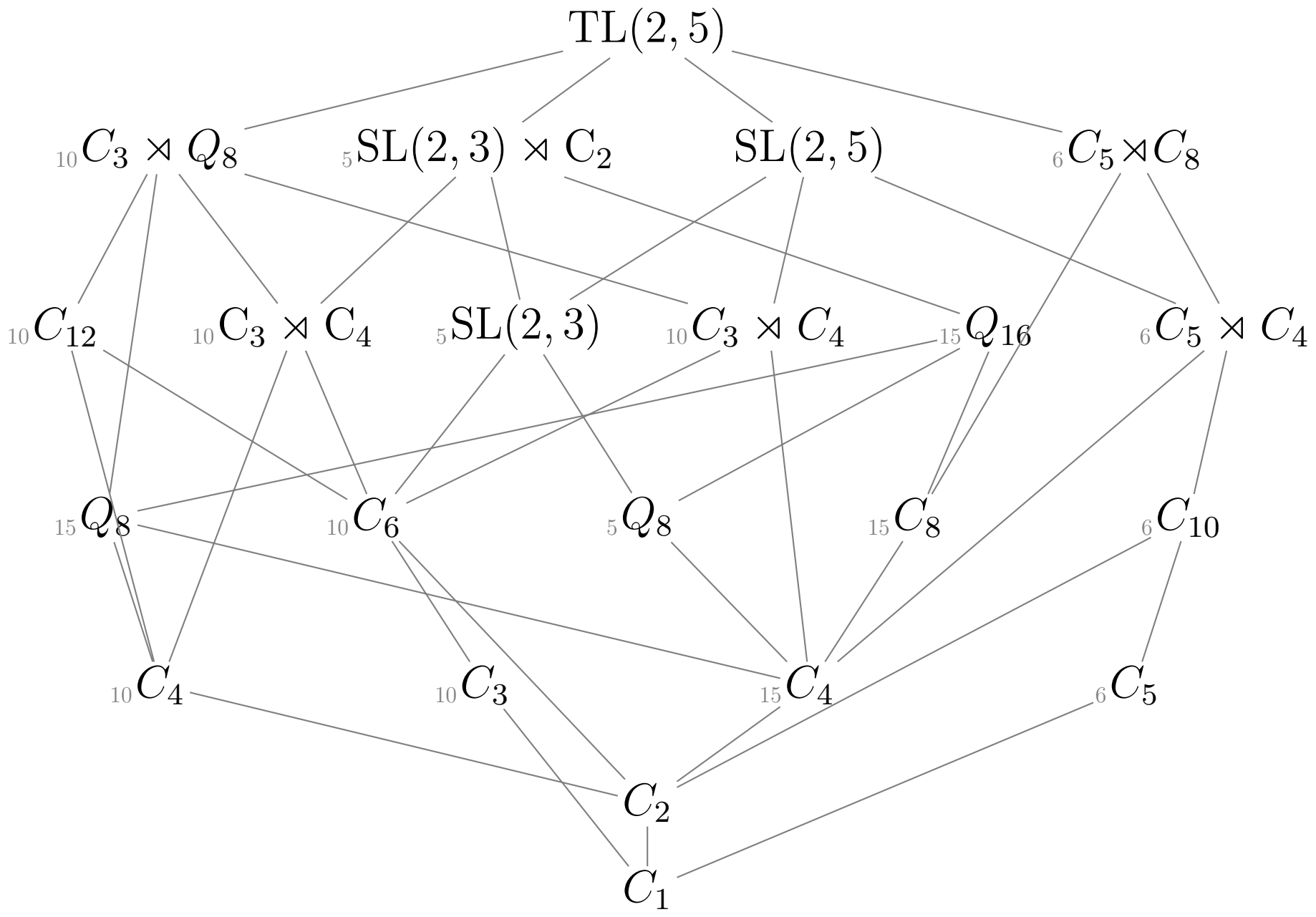
\begin{figure}[H]
  \centering
\begin{tikzpicture}[scale=1.0,sgplattice]
  \node at (5.5,0) (1) {\gn{C1}{C_1}};
  \node at (5.5,0.803) (2) {\gn{C2}{C_2}};
  \node at (4,1.89) (3) {\gn{C3}{C_3}};
  \node at (10,1.89) (4) {\gn{C5}{C_5}};
  \node at (1,1.89) (5) {\gn{C4}{C_4}};
  \node at (7,1.89) (6) {\gn{C4}{C_4}};
  \node at (3,3.44) (7) {\gn{C6}{C_6}};
  \node at (10.5,3.44) (8) {\gn{C10}{C_{10}}};
  \node at (5.5,3.44) (9) {\gn{Q8}{Q_8}};
  \node at (0.5,3.44) (10) {\gn{Q8}{Q_8}};
  \node at (8,3.44) (11) {\gn{C8}{C_8}};
  \node at (2.25,5.2) (12) {\gn{Dic3}{\mathrm{C_{3}\rtimes C_{4}}}};
  \node at (6.62,5.2) (13) {\gn{Dic3}{C_{3}\rtimes C_{4}}};
  \node at (0.125,5.2) (14) {\gn{C12}{C_{12}}};
  \node at (10.9,5.2) (15) {\gn{Dic5}{C_{5}\rtimes C_{4}}};
  \node at (8.75,5.2) (16) {\gn{Q16}{Q_{16}}};
  \node at (4.38,5.2) (17) {\gn{SL(2,3)}{{\rm SL}(2,3)}};
  \node at (1,6.84) (18) {\gn{Dic6}{C_{3}\rtimes Q_{8}}};
  \node at (10,6.84) (19) {\gn{C5:C8}{C_5{\rtimes}C_8}};
  \node at (4,6.84) (20) {\gn{CSU(2,3)}{{\rm{SL}(2,3)\rtimes C_{2}}}};
  \node at (7,6.84) (21) {\gn{SL(2,5)}{{\rm SL}(2,5)}};
  \node at (5.5,7.94) (22) {\gn{CSU(2,5)}{{\rm TL}(2,5)}};
  \draw[lin] (1)--(2) (1)--(3) (1)--(4) (2)--(5) (2)--(6) (2)--(7) (3)--(7)
     (2)--(8) (4)--(8) (6)--(9) (5)--(10) (6)--(10) (6)--(11) (5)--(12)
     (7)--(12) (6)--(13) (7)--(13) (5)--(14) (7)--(14) (6)--(15) (8)--(15)
     (9)--(16) (10)--(16) (11)--(16) (7)--(17) (9)--(17) (12)--(18) (13)--(18)
     (14)--(18) (10)--(18) (15)--(19) (11)--(19) (12)--(20) (16)--(20) (17)--(20)
     (13)--(21) (15)--(21) (17)--(21) (18)--(22) (19)--(22) (20)--(22) (21)--(22);
  \node[cnj=3] {10};
  \node[cnj=4] {6};
  \node[cnj=5] {10};
  \node[cnj=6] {15};
  \node[cnj=7] {10};
  \node[cnj=8] {6};
  \node[cnj=9] {5};
  \node[cnj=10] {15};
  \node[cnj=11] {15};
  \node[cnj=12] {10};
  \node[cnj=13] {10};
  \node[cnj=14] {10};
  \node[cnj=15] {6};
  \node[cnj=16] {15};
  \node[cnj=17] {5};
  \node[cnj=18] {10};
  \node[cnj=19] {6};
  \node[cnj=20] {5};
\end{tikzpicture}
 \caption{Conjugacy classes of subgroups of $\TL(2,5)$.}
  \label{fig:tl25}
\end{figure}

\begin{table}[H]
\centering
{\fontsize{10pt}{12pt}\selectfont 
\begin{tabular}{|c||c|c|c|c|c|c|c|c|c|c|c|c|c|} \hline
& 
$C_1$ & $C_{2}$ & $C_{3}$ & 
$C_{4(10)}$ & $C_{4(15)}$ & $C_{5}$ & 
$C_{6}$  & $C_{8}$ & $Q_{8(5)}$  & $Q_{8(15)}$ &$C_{10}$ & $C_{3}\rtimes C_{4}$ & $\mathrm{C_{3}\rtimes C_{4}}$ \\ 
\hline
\hline

$V_{1}$ & 
$1$ & $1$ & $1$ & $0$ & $1$ & 
$1$ & $1$ & $0$ & $1$ & $0$ & $1$ & $0$ & $1$\\
\hline

$V_{4,1}$ & 
$4$ & $4$ & $2$ & $1$ & $2$ & 
$0$ & $2$ & $1$ & $1$ & $0$ & $0$ & $0$ & $1$\\
\hline 
$V_{4,2}$ & 
$4$ & $4$ & $2$ & $3$ & $2$ & 
$0$ & $2$ & $1$ & $1$ & $2$ & $0$ & $2$ & $1$\\
\hline

$V_{5,1}$ & 
$5$ & $5$ & $1$ & $3$ & $3$ & 
$1$ & $1$ & $1$ & $2$  & $2$ & $1$ & $1$ & $1$\\
\hline 

$V_{5,2}$ & 
$5$ & $5$ & $1$ & $2$ & $3$ & 
$1$ & $1$ & $2$ & $2$  & $1$ & $1$ & $0$ & $1$\\

\hline 

$V_{6}$ & 
$6$ & $6$ & $2$ & $3$ & $2$ & 
$2$ & $2$ & $1$ & $0$  & $1$ & $2$ & $1$ & $0$\\
\hline 

$V_{8,1}$ & 
$8$ & $0$ & $0$ & $0$ & $0$ & 
$0$ & $0$ & $0$ & $0$  & $0$ & $0$ & $0$ & $0$\\

\hline 

$V_{8,2}$ & 
$8$ & $0$ & $4$ & $0$ & $0$ & 
$0$ & $0$ & $0$ & $0$  & $0$ & $0$ & $0$ & $0$\\

\hline

$V_{12}$ & 
$12$ & $0$ & $4$ & $0$ & $0$ & 
$4$ & $0$ & $0$ & $0$  & $0$ & $0$ & $0$ & $0$\\

\hline 
\end{tabular}
}
\vskip 0.2cm 
{\fontsize{10pt}{12pt}\selectfont 
\begin{tabular}{|c||c|c|c|c|c|c|c|c|c|c|c|c|c|c|c|c|c|c|} \hline
 & 

$C_{12}$  & $Q_{16}$ & $C_{5}\rtimes C_{4}$ &  
$\mathrm{SL}(2,3)$ &  $C_{3}\rtimes Q_{8}$ & 
$C_{5}\rtimes C_{8}$ &$\mathrm{SL}(2,3)\rtimes C_{2}$&$\mathrm{SL}(2,5)$&$\mathrm{TL}(2,5)$
\\ \hline
\hline
$V_{1}$ & 
$0$ & $0$ & $1$ & $1$ & $0$ & 
$0$ & $0$ & $1$ & $0$ \\
\hline

$V_{4,1}$ & 
$1$ & $0$ & $0$ & $1$ & $0$ & 
$0$ & $0$ & $0$ & $0$ \\
\hline

$V_{4,2}$ & 
$1$ & $1$ & $0$ & $1$ & $1$ & 
$0$ & $1$ & $0$ & $0$ \\
\hline

$V_{5,1}$ & 
$1$ & $1$ & $1$ & $0$ & $1$ & 
$0$ & $0$ & $0$ & $0$ \\
\hline 

$V_{5,2}$ & 
$0$ & $1$ & $1$ & $0$ & $0$ & 
$1$ & $0$ & $0$ & $0$ \\
\hline

$V_{6}$ & 
$1$ & $0$ & $0$ & $0$ & $0$ & 
$0$ & $0$ & $0$ & $0$ \\
\hline

$V_{8,1}$ & 
$0$ & $0$ & $0$ & $0$ & $0$ & 
$0$ & $0$ & $0$ & $0$ \\

\hline 

$V_{8,2}$ & 
$0$ & $0$ & $0$ & $0$ & $0$ & 
$0$ & $0$ & $0$ & $0$ \\

\hline 

$V_{12}$ & 
$0$ & $0$ & $0$ & $0$ & $0$ & 
$0$ & $0$ & $0$ & $0$ \\
\hline

\end{tabular}
}
\vskip 0.2cm
\caption{Fixed point dimension table for the non-trivial real irreducible representations of $\TL(2,5)$. A representative of the second $C_3\rtimes C_4$ is a subgroup of $\SL(2,5)$.}\label{table:dims_tl25}
\end{table}

\begin{figure}[H]
  \centering
   \begin{tikzpicture}[scale=1.0,sgplattice]
  \node at (5.5,0) (1) {\gn{C1}{C_1}};
  \node at (8,0.803) (2) {\gn{C2}{C_2}};
  \node at (3,0.803) (3) {\gn{C2}{C_2}};
  \node at (8,2.14) (4) {\gn{C3}{C_3}};
  \node at (5.5,2.14) (5) {\gn{C5}{C_5}};
  \node at (3,2.14) (6) {\gn{C2^2}{C_2^2}};
  \node at (10.5,2.14) (7) {\gn{C2^2}{C_2^2}};
  \node at (0.5,2.14) (8) {\gn{C4}{C_4}};
  \node at (10.9,3.89) (9) {\gn{C6}{C_6}};
  \node at (8.75,3.89) (10) {\gn{S3}{D_6}};
  \node at (4.38,3.89) (11) {\gn{S3}{S_3}};
  \node at (2.25,3.89) (12) {\gn{D5}{D_{10}}};
  \node at (0.125,3.89) (13) {\gn{D4}{D_8}};
  \node at (6.62,3.89) (14) {\gn{A4}{A_4}};
  \node at (10,5.54) (15) {\gn{D6}{D_{12}}};
  \node at (1,5.54) (16) {\gn{F5}{C_5\rtimes C_4}};
  \node at (7,5.54) (17) {\gn{S4}{S_4}};
  \node at (4,5.54) (18) {\gn{A5}{A_5}};
  \node at (5.5,6.62) (19) {\gn{S5}{S_5}};
  \draw[lin] (1)--(2) (1)--(3) (1)--(4) (1)--(5) (3)--(6) (2)--(7) (3)--(7)
     (3)--(8) (2)--(9) (4)--(9) (2)--(10) (4)--(10) (3)--(11) (4)--(11)
     (3)--(12) (5)--(12) (6)--(13) (7)--(13) (8)--(13) (4)--(14) (6)--(14)
     (7)--(15) (9)--(15) (10)--(15) (11)--(15) (12)--(16) (8)--(16) (13)--(17)
     (14)--(17) (10)--(17) (12)--(18) (14)--(18) (11)--(18) (15)--(19) (16)--(19)
     (17)--(19) (18)--(19);
  \node[cnj=2] {10};
  \node[cnj=3] {15};
  \node[cnj=4] {10};
  \node[cnj=5] {6};
  \node[cnj=6] {5};
  \node[cnj=7] {15};
  \node[cnj=8] {15};
  \node[cnj=9] {10};
  \node[cnj=10] {10};
  \node[cnj=11] {10};
  \node[cnj=12] {6};
  \node[cnj=13] {15};
  \node[cnj=14] {5};
  \node[cnj=15] {10};
  \node[cnj=16] {6};
  \node[cnj=17] {5};
\end{tikzpicture}
 \caption{Conjugacy classes of subgroups of $S_5$.}
  \label{fig:s5}
\end{figure}
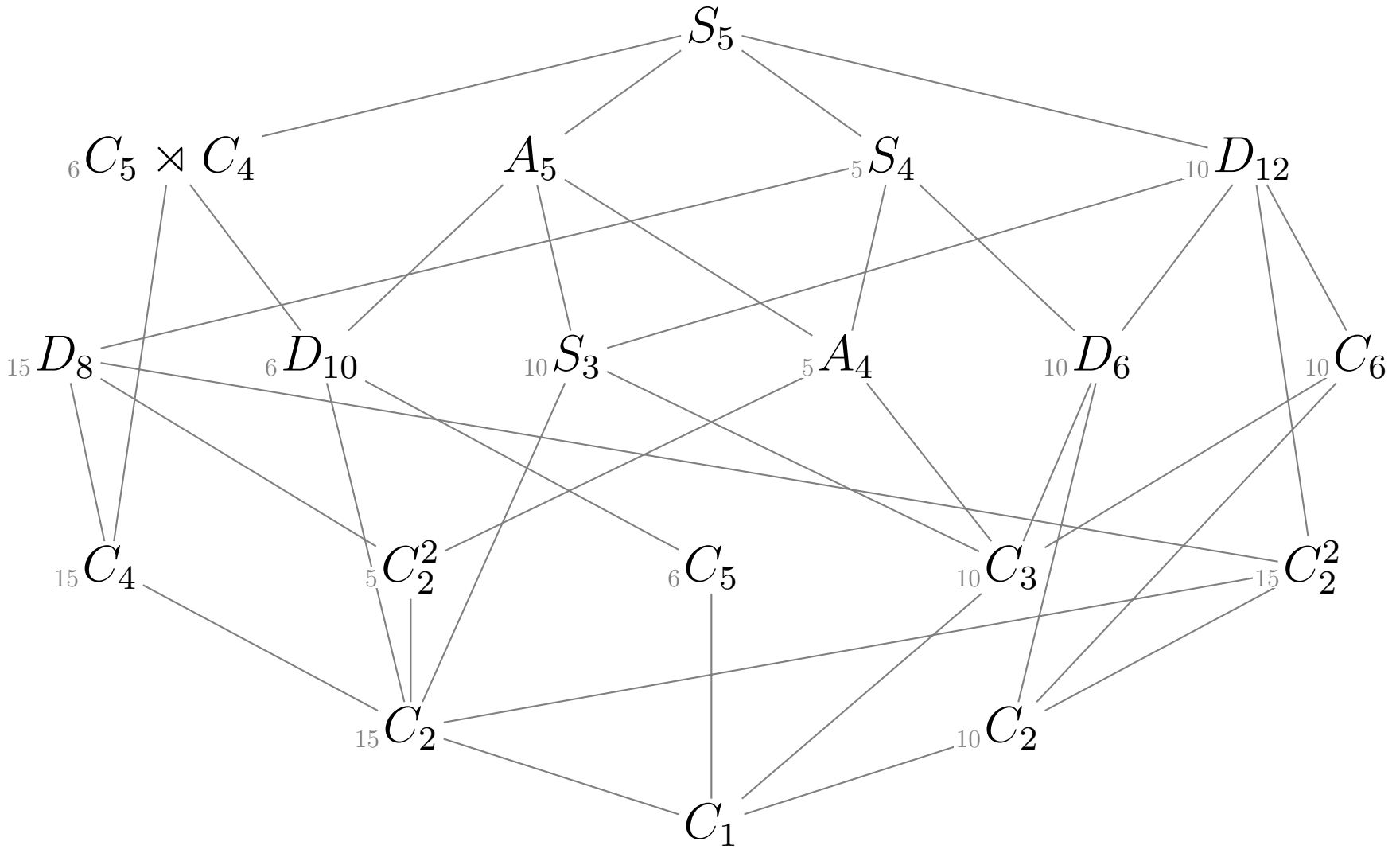

\begin{table}[H]
\centering
{\fontsize{10pt}{12pt}\selectfont 
\begin{tabular}{|c||c|c|c|c|c|c|c|c|c|c|c|c|c|c|c|c|c|c|c|} \hline
& 
$C_1$ & $(C_{2})_{10}$  &$(C_{2})_{15}$ & 
$C_{3}$  &$(C_{2}^2)_5$ & $C_{4}$ & $(C_{2}^2)_{15}$ &$C_{5}$& 

 $D_6$ & $S_3$ 
 \\ 
\hline
\hline

$V_{1}$ & 
$1$ & $0$ & $1$ & $1$ & $1$ & 
$0$ & $0$  & $1$ & 
 $0$ &$1$  \\
\hline 

$V_{4,1}$ & 
$4$ & $1$ & $2$ & $2$ & $1$ & 
$1$ & $0$ &$0$ & 
 $0$ &$1$  \\
\hline

$V_{4,2}$ & 
$4$ & $3$ & $2$ & $2$ & $1$ & 
$1$ & $2$ & $0$ & 
 $2$ &$1$  \\
\hline 

$V_{5,1}$ & 
$5$ &$3$& $3$ & $1$ & $2$ & $1$ & 
$2$ & $1$ & 
 $1$ &$1$   \\
\hline 

$V_{5,2}$ & 
$5$ & $2$ & $3$ & $1$ & $2$ & 
$2$ & $1$  & $1$ & 
 $0$ &$1$   \\

\hline

$V_{6}$ & 
$6$ & $3$ & $2$ & $2$ & $0$ & 
$1$ & $1$ & $2$ & 
 $1$ &$0$  \\
\hline

\end{tabular}
}
\vskip 0.2cm 
{\fontsize{10pt}{12pt}\selectfont 
\begin{tabular}{|c||c|c|c|c|c|c|c|c|c|c|c|c|c|c|c|c|c|c|} \hline
& 
$C_6$ & $D_8$  &$D_{10}$ & 
$A_4$  &$D_{12}$ & $C_5\rtimes C_4$ & $S_4$ &$A_5$& 

 $S_5$
 \\ 
\hline
\hline

$V_{1}$ & 
$0$ & $0$ & $1$ & $1$ & $0$ & 
$0$ & $0$  & $1$ & 
 $0$   \\
\hline 

$V_{4,1}$ & 
$1$ & $0$ & $0$ & $1$ & $0$ & 
$0$ & $0$ &$0$ & 
 $0$   \\
\hline

$V_{4,2}$ & 
$1$ & $1$ & $0$ & $1$ & $1$ & 
$0$ & $1$ & $0$ & 
 $0$   \\
\hline 
$V_{5,1}$ & 
$1$ &$1$& $1$ & $0$ & $1$ & $0$ & 
$0$ & $0$ & 
 $0$    \\
\hline 
$V_{5,2}$ & 
$0$ & $1$ & $1$ & $0$ & $0$ & 
$1$ & $0$  & $0$ & 
 $0$    \\

\hline

$V_{6}$ & 
$1$ & $0$ & $0$ & $0$ & $0$ & 
$0$ & $0$ & $0$ & 
 $0$   \\
\hline

\end{tabular}
}

\caption{Fixed point dimension table for the non-trivial real irreducible representations of $S_5$.}\label{table:dims_s5}
\end{table} 

In light of \Cref{corollary:fixed_dims2}, we list the real irreducible $G$-representations that arise from irreducible complex representations with Frobenius–Schur indicator equal to $1$. These representations are precisely those for which the dimension should be doubled in the second case of \Cref{corollary:fixed_dims2}. This avoids potential sources of confusion in the subsequent proofs. The lists are as follows:
\begin{itemize}
    \item $G = \SL(2,5)$: $V_3, V_{4,1}$, and $V_5$,
    \item $G = \TL(2,5)$: $V_1, V_{4,1}, V_{4,2}, V_{5,1}, V_{5,2}$, and $V_6$,
    \item $G = S_5$: $V_1, V_{4,1}, V_{4,2}, V_{5,1}, V_{5,2}$, and $V_6$.
\end{itemize}

\section{Exclusion techniques}\label{section:1}
We adapt the already known exclusion techniques for one fixed point actions to the case of direct products $G = H \times C_k$. The main idea is to use only the representation-theoretical information arising from $H$ to infer the nonexistence of the actions in question for $G$. Below, we recall the techniques mentioned above and show how they can be adapted to the direct products of interest.

Suppose that $G$ acts on a manifold $M$ with a fixed point $x$. Note that the tangent space $T_xM$ inherits an action of $G$ from $M$ via the differentials of the diffeomorphisms. In this way, $T_xM$ becomes a $G$-representation whose dimension equals that of $M$. The common starting point is the observation that the local behaviour around $x$ can be studied through the action of $G$ on the tangent space $T_xM$. This follows from the \emph{Slice Theorem}:

\begin{theorem}\emph{\cite[p. 308]{SliceThm}}\label{theorem:slice}
   There exists a $G$-invariant neighborhood $U$ of $x$ such that there exists a $G$-equivariant diffeomorphism $f\colon U\rightarrow T_x M$.
\end{theorem}

\noindent The following theorem provides general conditions for the nonexistence of one–fixed–point actions by means of intersection theory. In subsequent statements, $A_x$ denotes the connected component of a manifold $A$ containing $x\in A$.

\begin{theorem}\label{theorem:pdis}
Let $R$ be either the ring of integers $\mathbb{Z}$ or the finite field $\mathbb{Z}_{p}$, where $p$ is an arbitrary prime number. Suppose that a compact Lie group $G$ acts on a manifold $M$ with an isolated fixed point $x$ and tangential $G$-representation $T_x M \cong V$. Assume that there exist subgroups $H$, $H_1$, and $H_2$ of $G$ such that 
$H \leq H_1 \cap H_2$, the subgroups $H_1$ and $H_2$ generate $G$, the manifolds $M_x^{H}$, $M_x^{H_1}$, and $M_x^{H_2}$ are $R$-orientable, 
$\dim V^{H_1}, \dim V^{H_2} > 0$, 
and 
$\dim V^{H_1} + \dim V^{H_2} = \dim V^{H}$.
Denote by $i_j \colon M_x^{H_j} \rightarrow M_x^{H}$, for $j=1,2$, the inclusions of submanifolds, and let $\mathrm{PD}$ denote the Poincaré duality functor for $M_x^{H}$. Under these assumptions, if for some $R$ we have
$
\mathrm{PD}(i_{1\ast}[M_x^{H_1}]) \cup \mathrm{PD}(i_{2\ast}[M_x^{H_2}]) \neq \pm 1
$ in $H^{n}(M_x^{H}; R)$, $n=\dim V^H$, then $M^{G} \neq \{ \pt\}$.

\end{theorem}

\begin{proof}
    By interpreting the cup product in terms of the intersection number for submanifolds of complementary dimensions in $M_x^H$, we conclude that $M_x^{H_{1}} \cap M_x^{H_{2}}$ is not a single point. Since $x$ already belongs to $M_x^{H_{1}} \cap M_x^{H_{2}}$, there must exist another point in $M_x^{H_{1}} \cap M_x^{H_{2}} \subseteq M^{G}$.

\end{proof}

\noindent The theorem above will be applied by us primarily in the following cases:
\begin{enumerate}
    \item $H$ is a $2$-group.
    \item $H$ is a $p$-group for some odd prime $p$, and both $M^{H_1}$ and $M^{H_2}$ are orientable. The latter condition holds whenever the orders $|H_1|$ and $|H_2|$ are odd.
\end{enumerate}

\begin{theorem}\label{theorem:equaldims_Kamil}
    Let a compact Lie group $G$ act on a manifold $M$ with an isolated fixed point $x$ and tangent $G$-representation $T_x M\cong V$, such that there exist subgroups $H$, $H_{1}$ and $H_{2}$ of $G$ with $H \leq H_{1} \cap H_{2}$, $H_{1}$ and $H_{2}$ generate $G$, and $\dim V^{H_1} =\dim V^{H}$. Suppose further that at least one of the following conditions holds:
\begin{enumerate}
\item $M^{H}$ is connected and $H_{2}$ cannot act on $M$ with exactly one fixed point such that $T_{x}M$ is isomorphic to $\mathrm{Res}^G_{H_2}V$ as an $H_2$-representation,
\item $H\trianglelefteq H_{2}$ and $H_{2}/H$ cannot act on $M_x^H$ with exactly one fixed point such that $T_{x}M^H$ is isomorphic to $\mathrm{Res}^G_{H_2}V^H$ as an $H_{2}/H$-representation.
\end{enumerate}
Then the fixed point set $M^{G}$ contains more than one point.
\end{theorem}

\begin{proof}
    Suppose first that $M^{H}$ is connected and that $H_{2}$ cannot act on $M$ with exactly one fixed point such that $T_{x}M$ is isomorphic to $\mathrm{Res}^G_{H_2}V$ as an $H_2$ representation. Then, by the Principal Orbit Theorem, we have $M^{H} = M^{H_1}$, and $M^{H_2} \subset M^{H_1}$. Since $M^{G} = M^{H_1} \cap M^{H_2} = M^{H_2} \neq \{\pt\}$, the conclusion follows.

Now assume that $H\trianglelefteq H_{2}$ and that $H_{2}/H$ cannot act on $M_x^H$ with exactly one fixed point and with $T_{x}M^H$ isomorphic to $\mathrm{Res}^G_{H_2}V^H$ as an $H_{2}/H$ representation. Again, by the Principal Orbit Theorem, we know that $M_{x}^{H} = M_{x}^{H_1}$. For the contrary, suppose that $M_x^G$ is a single point. Then $M_x^G = M_x^{H_1}\cap M_x^{H_2}$, since $H_1$ and $H_2$ generate $G$. On the other hand, $M_x^{H_1}\cap M_x^{H_2} = M_x^{H}\cap M_x^{H_2} = M_x^{H_2} = (M_x^H)^{H_2/H}\neq\{\pt\}$.
\end{proof}

\begin{lemma}\label{lemma:8cond}
Let $H$ be a subgroup of $G$ such that each element of $H$ is of prime-power order. Suppose that $G$ acts on a $\mathbb{Z}$-homology sphere $\Sigma$ with a fixed point $x$ and tangent $G$-representation $T_{x}\Sigma \cong V$ such that $V^{H} = 0$. Then $\Sigma^{H}$ and $\Sigma^{G}$ are discrete sets. Suppose moreover that for every element $h \in H$ of $2$-power order, the order of $h$ is less than or equal to $4$ (cf. the $8$-condition in \emph{\cite{pawalowski_8_cond}}) then $\mathrm{Res}^{G}_{H}(T_{y}\Sigma) \cong \mathrm{Res}^{G}_{H}(V)$ for any $y \in \Sigma^{H}$.
\end{lemma}

\begin{proof}
 $\Sigma^{H}$ is a discrete set and  with isomorphic tangential $H$-representations by theorem \cite[Proposition 2.9]{MorimotoTamura}. The obvious inclusion $\Sigma^{G} \subseteq \Sigma^{H}$ concludes the proof.
\end{proof}

\noindent There is a similar statement that holds under slightly different assumptions.

\begin{lemma}\emph{(cf. \cite[Proposition 2.5]{morimoto_appendix_mizerka})}\label{lemma:equal_dims}
   Let $G$ act on a homology sphere $\Sigma$ with a fixed point $x$ and representation $T_x\Sigma$ isomorphic to $V$, and assume that there exist  subgroups $H,H_1$ and $H_2$ such that:
     \begin{enumerate}
    \item $H$ is an extension of a cyclic group by a $p$-group and each element of $H$ is of prime-power order,
    \item $H\leq H_{1}\cap H_{2}$ and $H_1$ and $H_2$ generate $G$,
    \item $\dim V^H =0$,
    \item $H_{i}$ cannot act on $\Sigma$ with exactly one fixed point and the representation at the fixed point isomorphic to $\mathrm{Res}^{G}_{H_{i}}V$ for $i=1,2$.
\end{enumerate}
     Then $\Sigma^G = S^{0}$ .

\end{lemma}
\begin{proof}
    By \Cref{lemma:8cond} $\Sigma^{G}$ is a discrete set, and by \cite[PROPOSITION 1]{Oliver1975} the Euler characteristic of $\Sigma^{H}\setminus\{x\}$ is equal to $1$, hence $\Sigma^{H} = S^{0}$. Since $\Sigma^{H_{i}} \neq \{\pt\}$, we have $\Sigma^{H_{i}}=\Sigma^{H}$ for $i=1,2$, and the result follows from the equality $\Sigma^{G} = \Sigma^{H_{1}}\cap \Sigma^{H_{2}} = \Sigma^{H} \cap \Sigma^{H} = \Sigma^{H}$.

\end{proof}

\noindent The following theorem is due to Conner and Floyd.
\begin{theorem}\emph{\cite[(25.1) Theorem]{ConnerFloyd}}\label{theorem:conner_floyd}
Suppose $C_2$ acts on a closed connected manifold $M$. Then, if $M^{C_2}$ is a non-empty discrete set, $M^{C_2}$ consists of an even number of points.

\end{theorem}

\noindent The following lemma is another very useful constraint.

\begin{lemma}\emph{\cite[Proposition 2.7]{morimoto_appendix_mizerka}}\label{lemma:one_dim}
Let $G$ act on a connected closed manifold $M$. If for some $x \in M^{G} $ 
the tangential representation $T_{x}M$ 
contains a non-trivial irreducible 
representation of dimension one,
then $M^{G}$ is not a singleton.

\end{lemma}

\begin{proof}
    Consider a non-trivial irreducible component of dimension one
    of $T_{x}M$ which can be thought of as a representation  
    $\varphi\colon G \rightarrow O(1) \cong C_{2}$.
    Then $\mathrm{Ker}\varphi$ 
    is a normal subgroup of index $2$,
    and we can write $M^{G}$ as 
    $(M^{\mathrm{Ker}\varphi})^{C_{2}}$. Suppose for the contrary that $M^G=\{x\}$.
    Using \Cref{theorem:slice} we obtain $\mathrm{dim}(M_x^{\mathrm{Ker\varphi}})=\dim(T_xM)^{\mathrm{Ker\varphi}}\geq 1>0$,
    and the result follows by \Cref{theorem:conner_floyd} by applying it to the action of $G/\mathrm{Ker\varphi}\cong C_2$ on $M_x^{\mathrm{Ker\varphi}}$.

\end{proof}

\noindent In multiple cases we shall infer exclusions of one fixed point actions by studying the behavior of subgroups whose index is a power of $2$. To that end we prove the following below.

\begin{theorem}\label{lemma:index_two_excl2}
    Let a compact Lie group $G$ 
    act on an $n$-dimensional $\mathbb{Z}_{2}$-homology sphere $\Sigma$, $n\neq 2$,
    with an isolated fixed point $x$
    and tangential $G$-representation $T_x\Sigma \cong V$. 
    Suppose further that there exist closed subgroups 
    $H_1, H_2$, and $H$ of $G$ such that:
    \begin{itemize}
    \item $\langle H_1, H_2\rangle =H$, 
    \item $H_i\trianglelefteq G$ for $i=1$ and $i=2$ and $H \trianglelefteq G$ and $|G/H| = 2^{k}$ for some nonnegative $k$,
    \item $\dim V^{H_1}+\dim V^{H_2}=n$ and either $\dim V^{H_1}, \dim V^{H_2}>0$ or if $\dim V^{H_i}=0$ then $\Sigma^{H_i}$ contains an even number of points.
\end{itemize}
  Then $\Sigma^G$ is not a singleton.

\end{theorem}

 \begin{proof}
    Note first that in the case $\dim V^{H_i}=0$, $\Sigma^{H_i}$ is finite and contains an even number of points by the third assumption. Thus, the $2$-group $G/H$ acts on $\Sigma^H=\Sigma^{H_1}\cap\Sigma^{H_2}=\Sigma^{H_i}$ with one fixed point. A contradiction.
    
    We can assume therefore $\dim V^{H_1}, \dim V^{H_2}>0$. Let
    $$A =\Sigma_{x}^{H_1} \cap \Sigma_{x}^{H_2} = \bigcup_{i \in I}\{x_{i}\} \cup \bigcup_{j \in J}Y_{j},$$ 
    be the decomposition of this part of the fixed point set $\Sigma^{H}$ which is contained in both $\Sigma_{x}^{H_1} $ and $\Sigma_{x}^{H_2}$, 
    where by $x_{i}$ we denote the $H$-isolated fixed points, 
    and by $Y_{j}$ the components of strictly positive dimension, 
    and both indexing sets $I$ and $J$ are finite.

By assumption that $\Sigma$ is a $\mathbb{Z}_{2}$-homology sphere we have 
$
\mathrm{PD}(i^*(\Sigma_x^{H_{1}})) \cup \mathrm{PD}(i^*(\Sigma_x^{H_{2}})) = 0 \mod 2.
$
Now, by removing the $H$-tubular neighborhoods of all isolated $H$-fixed points which are in $A$,
we obtain a manifold with boundary $(X, \partial X)$ such that
$$
H^i(X, \partial X; \mathbb{Z}_2) \cong 
\begin{cases}
\mathbb{Z}_2 & i = n, \\
\mathbb{Z}_2^{k-1} & i = 1, \\
0 & \text{otherwise.}
\end{cases}
$$
where $k$ is the number of isolated fixed points.
By the relative Thom transversality theorem, we may perturb
$\Sigma_x^{H_1}\cap X$ slightly, relative to a neighbourhood of $\partial X$, so that it becomes transverse to
$\Sigma_x^{H_2}\cap X$. The perturbation is chosen to be supported only in small neighbourhoods of the positive-dimensional components of
$$
        (\Sigma_x^{H_1}\cap X)\cap(\Sigma_x^{H_2}\cap X),
$$
and it is kept fixed near the isolated intersection points. In particular, the perturbation does not change the local contributions coming from the isolated points removed from $\Sigma$.

From the ring structure
$
H^*(X, \partial X; \mathbb{Z}_2)
$
we obtain
$
\mathrm{PD}(i^*(\Sigma_x^{H_1} \cap X)) \cup \mathrm{PD}(i^*(\Sigma_x^{H_2} \cap X)) = 0 \mod 2.
$
By interpreting the cup product as the mod $2$ intersection number, the total intersection number of $\Sigma_x^{H_1}$ and $\Sigma_x^{H_2}$ in $\Sigma$ is zero modulo $2$. After the relative perturbation, this total intersection number is the sum, modulo $2$, of the contributions from the isolated points of $\Sigma^H$ and the transverse intersection points inside $X$ arising from the positive-dimensional part of the intersection. The latter contribution is even, since it is precisely the relative mod $2$ intersection number in $X$. Hence the number of isolated points of $A$ is even. Moreover, this set of isolated points is $G/H$-invariant and contains the point $x$ and thus, it is nonempty. 
To finish the proof, it remains to observe that if a $2$-group $G/H$ acts on a finite set of even cardinality, then the number of fixed points is also even, because all non-fixed orbits have even cardinality. Since $x$ is fixed by $G/H$, the fixed point set is non-empty, and therefore it contains at least two points. Hence
$$
        \Sigma^G\neq \{\mathrm{pt}\}.
$$

     \end{proof}

  \begin{corollary}\label{corolarry:one_dims_struct}
      Let $H\trianglelefteq G$, $G/H\cong C_2^{n}$, and suppose that $G \times C_k$ acts on a $\mathbb{Z}_{2}$-homology sphere $\Sigma$ with $G \times C_k$-representation at the isolated fixed point $x$ isomorphic to $\bigoplus_{i=1}^{n}(V_{1,i}\otimes \mathbb{C}_{k})^{m_{i}} \oplus \mathbb{C}_{k}^{m} \oplus W$,
      where $V_{1,i}$ is a $1$-dimensional representation with non-trivial action of the $i$-th factor in $C_2^{n}$, and $W$ is an arbitrary $G$-representation with the trivial $C_{k}$-action. Then $\Sigma^{G\times C_{k}} \neq \pt$.

  \end{corollary}
  \begin{proof}
      Suppose $\Sigma^{G\times C_{k}} =\pt$. Obviously, $W$ cannot contain the trivial factor. If $W$ contains a factor of the form $V_{1,i}$, then the result follows from \Cref{lemma:one_dim}. Otherwise it is enough to show that $\dim W^{H} = \dim W ^{G} = 0$, and apply \Cref{lemma:index_two_excl2} to the subgroups $C_{k},H <  H\times C_{k} \triangleleft G\times C_{k}$. 

Suppose now that
$$
        \dim W^H>\dim W^G=0.
$$
Then $W^H$ is a non-zero real representation of the quotient group
$$
        G/H\cong C_2^n.
$$
Since $G/H \cong C_2^{n}$, every irreducible real $G/H$-representation is
one-dimensional. Since $\dim (W^H)^{G/H}=\dim W^G=0$, the representation $W^H$
contains a non-trivial one-dimensional irreducible $G/H$-summand.
  \end{proof}

\section{Proofs of the main theorems}\label{section:3}

\subsection{The case of odd $k>2$}\label{subsection:3_1}

\begin{lemma}\label{lemma:general_form_sl25}
The following representation $W$ of $\SL(2,5)\times C_k$:

     $$   V_{4,1}^{\oplus k_{1}}   \oplus V_{4,2}^{\oplus k_{2}} \oplus V_{5}^{\oplus k_{3}}   \oplus V_{8}^{\oplus k_{4}}   \oplus V_{12}^{\oplus k_{5}} \oplus \mathbb{C}_k^{\oplus k_{6}}\oplus (V_{4,2}\otimes_{\mathbb{C}} \mathbb{C}_k)^{\oplus k_{7}}  $$
     cannot occur as the tangential
representation at the fixed point of a one fixed point action
on a homology sphere, where $k_{i}$
is an arbitrary non-negative integer.

\end{lemma}

\begin{proof}
It follows from \Cref{fig:sl25} that 
        $ C_{4}= ((C_{3}\rtimes C_{4}) \times C_{k}) \cap Q_{8} $ 
        and $\langle  (C_{3}\rtimes C_{4}) \times C_{k} , Q_{8}  \rangle = \SL(2,5)\times C_k$. 
        Applying \Cref{table:dims_sl25} and \Cref{corollary:fixed_dims2}, we can compute the dimensions of the fixed point sets:
        $$
        \mathrm{dim}W^{ C_{4}}=  2k_{1} + 3k_{3}+2k_6,$$
         $$ \mathrm{dim}W^{(C_{3}\rtimes C_{4}) \times C_{k}}= k_{1} + k_{3}, $$
        $$ \mathrm{dim}W^{ Q_{8} } =k_{1} + 2k_{3}+2k_6, $$
       so $\mathrm{dim}W^{ C_{4}} =\mathrm{dim}W^{ Q_{8}} +\mathrm{dim}W^{(C_{3}\rtimes C_{4}) \times C_{k}}$.

       We split the proof into three cases.

       First assume that $k_{1}=k_{3}=k_{6}=0$. Then $\mathrm{dim}W^{C_{4}}=0$. We apply \Cref{lemma:equal_dims} to the subgroups
       $$
       H=C_4,\quad H_1=Q_8,\quad H_2=(C_3\rtimes C_4)\times C_k.
       $$
       The subgroup $C_4$ is a $2$-group, hence it is an extension of a cyclic group by a $2$-group, and all its elements have prime-power order. Moreover, $C_4\leq Q_8\cap ((C_3\rtimes C_4)\times C_k)$, the subgroups $Q_8$ and $(C_3\rtimes C_4)\times C_k$ generate $\SL(2,5)\times C_k$, and $\mathrm{dim}W^{C_4}=0$. Finally, both $Q_8$ and $(C_3\rtimes C_4)\times C_k$ cannot act with exactly one fixed point in the required tangential representation type. For the latter subgroup this follows from the non-Oliver sequence
       $$
       C_3\trianglelefteq C_3\times C_k\trianglelefteq (C_3\rtimes C_4)\times C_k.
       $$
       Thus \Cref{lemma:equal_dims} gives a contradiction.

       Now assume that $k_{1}=k_{3}=0$ and $k_6>0$. Then
       $$
       \mathrm{dim}W^{C_4}=2k_6=\mathrm{dim}W^{Q_8}.
       $$
       We apply \Cref{theorem:equaldims_Kamil} to the subgroups
       $$
       H=C_4,\quad H_1=Q_8,\quad H_2=(C_3\rtimes C_4)\times C_k.
       $$
       The subgroup inclusions and the generation condition are as above, and the equality $\mathrm{dim}W^{H_1}=\mathrm{dim}W^H$ holds. Since $H_2=(C_3\rtimes C_4)\times C_k$ cannot act with exactly one fixed point in the required tangential representation type, \Cref{theorem:equaldims_Kamil} gives a contradiction.

       Finally assume that $k_1+k_3>0$. Then
       $$
       \mathrm{dim}W^{Q_8}>0
       $$
       and
       $$
       \mathrm{dim}W^{(C_3\rtimes C_4)\times C_k}>0.
       $$
       We apply \Cref{theorem:pdis} to the subgroups
       $$
       H=C_4,\quad H_1=Q_8,\quad H_2=(C_3\rtimes C_4)\times C_k.
       $$
       The subgroup inclusions and the generation condition are as above. Thus \Cref{theorem:pdis} and the equality
       $$
       \mathrm{dim}W^{ C_{4}} =\mathrm{dim}W^{ Q_{8}} +\mathrm{dim}W^{(C_{3}\rtimes C_{4}) \times C_{k}}
       $$
       gives the required complementary-dimension condition. Therefore \Cref{theorem:pdis} gives a contradiction.

       In all cases we obtain a contradiction, and hence $W$ cannot occur as the tangential representation at the fixed point of a one fixed point action on a homology sphere.

\end{proof}

\begin{theorem}\label{theorem:1main2}\emph{(cf. \Cref{theorem:1main})}
$\SL(2,5)\times C_k$ does not admit a one fixed point action on $S^n$ if $n\in\{8,9\}$.
\end{theorem}
\begin{proof}

   Suppose $n=8$. Then we have five possibilities: 
   $$
   V\in\{V_{4,2}\oplus \mathbb{C}_k^{\oplus 2},
   V_{4,2}\oplus (V_{4,2}\otimes_{\mathbb{C}}\mathbb{C}_k),
   (V_{4,2}\otimes_{\mathbb{C}}\mathbb{C}_k)^{\oplus 2},
   (V_{4,2}\otimes_{\mathbb{C}}\mathbb{C}_k)\oplus \mathbb{C}_k^{\oplus 2},
   V_8\otimes_{\mathbb{C}}\mathbb{C}_k\}.
   $$
   In the case $V=V_{4,2}\oplus \mathbb{C}_k^{\oplus 2}$, the assertion follows by \Cref{lemma:general_form_sl25}. Let $V\in\{V_{4,2}\oplus (V_{4,2}\otimes_{\mathbb{C}}\mathbb{C}_k),(V_{4,2}\otimes_{\mathbb{C}}\mathbb{C}_k)^{\oplus 2}\}$. 
   By \Cref{fig:sl25}, \Cref{table:dims_sl25}, and \Cref{corollary:fixed_dims2} we note that $\langle C_6,C_{10}\times C_k\rangle=\SL(2,5)\times C_k$, $C_2=C_6\cap(C_{10}\times C_k)$, and $\dim V^{C_2}=0$, and the assertion holds for these cases by 
   \Cref{lemma:equal_dims} and the fact that $C_{10}\times C_k$ is a non-Oliver group. Indeed, the following sequence of subgroups shows that $C_{10}\times C_k$ 
   is not Oliver: $C_2\trianglelefteq C_2\times C_k\trianglelefteq C_{10}\times C_k$.

   In the case $V=(V_{4,2}\otimes_{\mathbb{C}}\mathbb{C}_k)\oplus \mathbb{C}_k^{\oplus 2}$, we observe that $\dim V^{C_2}=\dim V^{C_6}$. We apply \Cref{theorem:equaldims_Kamil} to the subgroups
$$
H=C_2,\quad H_1=C_6,\quad H_2=C_{10}\times C_k.
$$
The subgroup $C_{10}\times C_k$ is not an Oliver group, as follows from the sequence
$$
C_2\trianglelefteq C_2\times C_k\trianglelefteq C_{10}\times C_k.
$$
Thus \Cref{theorem:equaldims_Kamil} gives a contradiction.
   

   Finally, consider the case $V=V_8\otimes_{\mathbb{C}}\mathbb{C}_k$. By \Cref{fig:sl25}, \Cref{table:dims_sl25}, and \Cref{corollary:fixed_dims2}, we have
   $$
   C_2=C_{10}\cap(C_6\times C_k),\quad
   \langle C_{10},C_6\times C_k\rangle=\SL(2,5)\times C_k,
   $$
   and
   $$
   \dim V^{C_2}=0.
   $$
   Hence the assertion follows from \Cref{lemma:equal_dims} applied to the subgroups
   $$
   H=C_2,\quad H_1=C_{10},\quad H_2=C_6\times C_k.
   $$
  
   Suppose $n=9$. Then we have three possibilities: 
   $$
   V\in\{V_3\oplus V_{4,2}\oplus\mathbb{C}_k,V_3\oplus (V_{4,2}\otimes_{\mathbb{C}}\mathbb{C}_k)\oplus\mathbb{C}_k,(V_{4,2}\otimes_{\mathbb{C}}\mathbb{C}_k)\oplus V_5\}.$$
   In the latter case, the assertion follows from \Cref{lemma:general_form_sl25}. In the first two cases, on the other hand, we note that $\dim V^{Q_8}=2$ and by the Smith Theory, combined with the classification of $2$-surfaces, we conclude that $\Sigma^{Q_8}\cong S^2$. On the other hand, we have $\dim V^{\SL(2,5)}=2$. Thus, by the Principal Orbit Theorem, it follows that $\Sigma^{\SL(2,5)}$ is an open and dense subset of $\Sigma^{Q_8}$. Being a closed submanifold of $\Sigma^{Q_8}$ at the same time, it is apparent that $\Sigma^{\SL(2,5)}\cong S^2$. Finally, $ \Sigma^{\SL(2,5) \times C_{k}}
 =(\Sigma^{\SL(2,5) })^{C_{k}} 
 \cong (S^{2})^{C_{k}} \neq \{\pt\}$.

\end{proof}

\begin{lemma}\label{lemma:structural_tl25}
Let $W$ be one of the following representations of $\TL(2,5)\times C_k$:
\begin{enumerate}
    \item $ V_{4,1}^{\oplus k_{1}}  \oplus V_{4,2}^{\oplus k_{2}} \oplus V_{5,1}^{\oplus k_{3}}   \oplus V_{5,2}^{\oplus k_{4}}   \oplus V_{8,1}^{\oplus k_{5}} \oplus V_{8,2}^{\oplus k_{6}}  \oplus V_{12}^{\oplus k_{7}}\oplus\mathbb{C}_k^{\oplus k_{8}} \oplus (V_{1}\otimes \mathbb{C}_k)^{\oplus k_{9}}\oplus (V_{4,1}\otimes \mathbb{C}_k)^{\oplus k_{10}} \oplus (V_{5,2}\otimes \mathbb{C}_k)^{\oplus k_{11}} \oplus (V_{8,1}\otimes_{\mathbb{C}}\mathbb{C}_k)^{k_{12}}\oplus(V_{8,2}\otimes_{\mathbb{C}}\mathbb{C}_k)^{k_{13}}\oplus (V_{12}\otimes_{\mathbb{C}}\mathbb{C}_k)^{k_{14}},  $
    \item  $V_{4,2}^{\oplus k_{1}} \oplus V_{5,1}^{\oplus k_{2}}   \oplus V_{5,2}^{\oplus k_{3}} \oplus V_{8,1}^{\oplus k_{4}}\oplus V_{8,2}^{\oplus k_{5}}   \oplus V_{12}^{\oplus k_{6}} \oplus\mathbb{C}_k^{\oplus k_{7}} \oplus (V_{1}\otimes \mathbb{C}_k)^{\oplus k_{8}}  \oplus (V_{4,2}\otimes \mathbb{C}_k)^{\oplus k_{9}} \oplus (V_{5,1}\otimes \mathbb{C}_k)^{\oplus k_{10}} \oplus(V_{8,1}\otimes_{\mathbb{C}}\mathbb{C}_k)^{k_{11}}\oplus(V_{8,2}\otimes_{\mathbb{C}}\mathbb{C}_k)^{k_{12}} \oplus (V_{12}\otimes_{\mathbb{C}}\mathbb{C}_k)^{k_{13}}.$
\end{enumerate}
Then $W$ cannot occur as the tangential
representation at the fixed point of a one fixed point action
on a homology sphere, where $k_{i}$
are arbitrary non-negative integers.

\end{lemma}

\begin{proof}
Assume, for a contradiction, that there exists a homology sphere $\Sigma$ with a
one-fixed-point action of $\TL(2,5)\times C_k$ and with the tangential representation at the
fixed point isomorphic to $W$.

We consider first case $(1)$. Recall from \Cref{fig:tl25} that
$$
        {Q_8}_{(15)}=Q_{16}\cap ((C_3\rtimes Q_8)\times C_k)
$$
and
$$
        \langle Q_{16},(C_3\rtimes Q_8)\times C_k\rangle
        =
        \TL(2,5)\times C_k.
$$
Put
$$
        H={Q_8}_{(15)},\qquad
        H_1=Q_{16},\qquad
        H_2=(C_3\rtimes Q_8)\times C_k.
$$
Using \Cref{table:dims_tl25} and \Cref{corollary:fixed_dims2}, we compute
$$
        \dim W^{{Q_8}_{(15)}}
        =
        2k_2+2k_3+k_4+2k_8+2k_{11},
$$
$$
        \dim W^{(C_3\rtimes Q_8)\times C_k}
        =
        k_2+k_3,
$$
and
$$
        \dim W^{Q_{16}}
        =
        k_2+k_3+k_4+2k_8+2k_{11}.
$$
Therefore
$$
        \dim W^{{Q_8}_{(15)}}
        =
        \dim W^{Q_{16}}
        +
        \dim W^{(C_3\rtimes Q_8)\times C_k}.
$$
We divide the proof into three cases.

First suppose that
$$
        \dim W^{{Q_8}_{(15)}}=0.
$$
Then also
$$
        \dim W^{Q_{16}}=
        \dim W^{(C_3\rtimes Q_8)\times C_k}=0.
$$
The subgroup $H={Q_8}_{(15)}$ is a $2$-group, hence it satisfies the group-theoretic
hypothesis of \Cref{lemma:equal_dims}. Moreover, $H\leq H_1\cap H_2$, the subgroups $H_1$ and
$H_2$ generate $\TL(2,5)\times C_k$, and $\dim W^H=0$. Finally, $Q_{16}$ is not an
Oliver group, and $(C_3\rtimes Q_8)\times C_k$ is not an Oliver group either. Indeed,
the latter follows from the sequence
$$
        C_3 \triangleleft C_3\times C_k \triangleleft (C_3\rtimes Q_8)\times C_k.
$$
Thus \Cref{lemma:equal_dims} gives a contradiction.

Now suppose that
$$
        \dim W^{{Q_8}_{(15)}}=\dim W^{Q_{16}}>0.
$$
Since ${Q_8}_{(15)}$ is a $2$-group and $\dim W^{{Q_8}_{(15)}}>0$, Smith theory implies
that the connected component $\Sigma^{{Q_8}_{(15)}}_x$ ($x$ denotes the single fixed point of $\TL(2,5)\times C_k$) is a positive-dimensional
$\mathbb Z_2$-homology sphere. In particular, it is connected. Hence the hypotheses of
\Cref{theorem:equaldims_Kamil} are satisfied for
$$
        H={Q_8}_{(15)},\qquad
        H_1=Q_{16},\qquad
        H_2=(C_3\rtimes Q_8)\times C_k.
$$
Indeed, $H\leq H_1\cap H_2$, the subgroups $H_1$ and $H_2$ generate
$\TL(2,5)\times C_k$, and
$$
        \dim W^H=\dim W^{H_1}.
$$
The subgroup $H_2$ cannot act with exactly one fixed point, because it is not an Oliver
group. Therefore \Cref{theorem:equaldims_Kamil} gives a contradiction.

It remains to consider the case
$$
        \dim W^{Q_{16}}>0
        \qquad\text{and}\qquad
        \dim W^{(C_3\rtimes Q_8)\times C_k}>0.
$$
Then the equality
$$
        \dim W^{{Q_8}_{(15)}}
        =
        \dim W^{Q_{16}}
        +
        \dim W^{(C_3\rtimes Q_8)\times C_k}
$$
gives the complementary-dimension condition required in \Cref{theorem:pdis}. Since
${Q_8}_{(15)}$ is a $2$-group, we apply \Cref{theorem:pdis} with coefficients in $\mathbb Z_2$.
All the remaining hypotheses have already been checked above: we have
$H\leq H_1\cap H_2$ and $\langle H_1,H_2\rangle=\TL(2,5)\times C_k$. Hence \Cref{theorem:pdis}
gives a contradiction.

We now consider case $(2)$. This time we use the subgroups
$$
        H=C_8,\qquad
        H_1=Q_{16},\qquad
        H_2=(C_5\rtimes C_8)\times C_k.
$$
Recall from \Cref{fig:tl25} that
$$
        C_8=Q_{16}\cap ((C_5\rtimes C_8)\times C_k)
$$
and
$$
        \langle Q_{16},(C_5\rtimes C_8)\times C_k\rangle
        =
        \TL(2,5)\times C_k.
$$
Using \Cref{table:dims_tl25} and \Cref{corollary:fixed_dims2}, we obtain
$$
        \dim W^{C_8}
        =
        k_1+k_2+2k_3+2k_7+2k_9+2k_{10},
$$
$$
        \dim W^{(C_5\rtimes C_8)\times C_k}
        =
        k_3,
$$
and
$$
        \dim W^{Q_{16}}
        =
        k_1+k_2+k_3+2k_7+2k_9+2k_{10}.
$$
Thus
$$
        \dim W^{C_8}
        =
        \dim W^{Q_{16}}
        +
        \dim W^{(C_5\rtimes C_8)\times C_k}.
$$
Again we split the proof into cases.

If
$$
        \dim W^{C_8}=0,
$$
then \Cref{lemma:equal_dims} applies to
$$
        H=C_8,\qquad
        H_1=Q_{16},\qquad
        H_2=(C_5\rtimes C_8)\times C_k.
$$
Indeed, $C_8$ is a $2$-group, $H\leq H_1\cap H_2$, the subgroups $H_1$ and $H_2$
generate $\TL(2,5)\times C_k$, and both $Q_{16}$ and $(C_5\rtimes C_8)\times C_k$
are non-Oliver groups. For the latter subgroup, one can use the sequence
$$
        C_5 \triangleleft C_5\times C_k \triangleleft (C_5\rtimes C_8)\times C_k.
$$
Hence \Cref{lemma:equal_dims} gives a contradiction.

If
$$
        \dim W^{C_8}=\dim W^{Q_{16}}>0,
$$
then
$$
        \dim W^{(C_5\rtimes C_8)\times C_k}=0.
$$
Since $C_8$ is a $2$-group and $\dim W^{C_8}>0$, Smith theory implies that
$\Sigma^{C_8}_x$ is a connected positive-dimensional $\mathbb Z_2$-homology sphere.
Therefore \Cref{theorem:equaldims_Kamil} applies to the subgroups
$$
        H=C_8,\qquad
        H_1=Q_{16},\qquad
        H_2=(C_5\rtimes C_8)\times C_k.
$$
Indeed, $H\leq H_1\cap H_2$, the subgroups $H_1$ and $H_2$ generate
$\TL(2,5)\times C_k$, and
$$
        \dim W^H=\dim W^{H_1}.
$$
Since $H_2$ is not an Oliver group, it cannot act with exactly one fixed point. Thus
\Cref{theorem:equaldims_Kamil} gives a contradiction.

Finally, suppose that
$$
        \dim W^{Q_{16}}>0
        \qquad\text{and}\qquad
        \dim W^{(C_5\rtimes C_8)\times C_k}>0.
$$
Then
$$
        \dim W^{C_8}
        =
        \dim W^{Q_{16}}
        +
        \dim W^{(C_5\rtimes C_8)\times C_k},
$$
so the complementary-dimension condition of \Cref{theorem:pdis} is satisfied. Since $C_8$ is
a $2$-group, we apply \Cref{theorem:pdis} with coefficients in $\mathbb Z_2$. This gives the
final contradiction.

In all cases we obtain a contradiction. Therefore $W$ cannot occur as the tangential
representation at the fixed point of a one-fixed-point action on a homology sphere.
\end{proof}

\begin{theorem}\label{theorem:2main2}\emph{(cf. \Cref{theorem:2main})}
$\TL(2,5)\times C_k$ does not admit a one fixed point action on a homology sphere $\Sigma$ of dimension $n$ where $n\in\{6,\ldots,13\}\cup\{15,16,17,21\}$.  
  
\end{theorem}

\begin{proof}
Suppose $\TL(2,5)\times C_k$ acts with exactly one fixed point $x$ on a homology sphere $\Sigma$ of dimension $n$ and let $V$ be the representation at $x$. Combining both cases of \Cref{lemma:structural_tl25} and \Cref{lemma:one_dim} we can conclude that $V$ must contain irreducible summands isomorphic to summands from each of the sets
$$
\mathcal{S}_1=\{V_{6},V_{4,2}\otimes\mathbb{C}_k,V_{5,1}\otimes\mathbb{C}_k,V_{6}\otimes \mathbb{C}_k\},
$$
and
$$
\mathcal{S}_2=\{V_{4,1},V_{6},V_{4,1}\otimes\mathbb{C}_k,V_{5,2}\otimes\mathbb{C}_k,V_{6}\otimes \mathbb{C}_k\}.
$$
Since the original action of $\TL(2,5)\times C_k$ is effective, the restriction of the tangential representation $V$ to $\TL(2,5)$ has to be effective. Indeed, suppose that some non-trivial element $g\in \TL(2,5)$ acts trivially on $T_x\Sigma$. By \Cref{theorem:slice}, $g$ acts trivially on a $\TL(2,5)\times C_k$-invariant neighbourhood of $x$. Since the principal orbit type is open and dense, and since the fixed point set of $g$ is closed, it follows that $g$ acts trivially on all of $\Sigma$. This contradicts the effectiveness of the original action. Hence $\mathrm{Res}^{\TL(2,5)\times C_k}_{\TL(2,5)}V$ is effective.

By the description of the irreducible $\TL(2,5)$-representations used in \Cref{table:dims_tl25}, all irreducible $\TL(2,5)\times C_k$-representations except the following ones have non-trivial kernel (we can read off the intersection of such kernel and $\TL(2,5)$ directly from \Cref{table:dims_tl25}):
$$
V_{8,i},\quad V_{8,j}\otimes_{\mathbb{C}}\mathbb{C}_k,\quad V_{12},\quad V_{12}\otimes_{\mathbb{C}}\mathbb{C}_k,\qquad i,j\in\{1,2\}.
$$
Therefore $V$ must contain a summand from the set 
$$
\mathcal{F}=\{V_{8,i},V_{8,j}\otimes_{\mathbb{C}}\mathbb{C}_k,V_{12},V_{12}\otimes_{\mathbb{C}}\mathbb{C}_k\},\quad i,j\in\{1,2\}.
$$
The above restrictions already exclude the cases $n\in\{6,\ldots,13\}\cup\{15,17\}$. We are left therefore with the cases $n=16$ and $n=21$. 

By \Cref{lemma:structural_tl25}, the $16$-dimensional representations $V$ which can appear 
as the tangential representation at the fixed point are of the forms:
$$
 V_6\oplus V_{8,i}\oplus\mathbb{C}_k,\quad V_6\oplus V_{8,i}\oplus(V_1\otimes\mathbb{C}_k),\quad
 V_6\oplus(V_{8,i}\otimes_{\mathbb{C}}\mathbb{C}_k)\oplus\mathbb{C}_k,\quad V_6\oplus(V_{8,i}\otimes_{\mathbb{C}}\mathbb{C}_k)\oplus(V_1\otimes\mathbb{C}_k),
$$
and the $21$-dimensional representations $V$, in turn, are of the forms:
$$
  V_6\oplus V_{8,i}\oplus V_{5,j}\oplus(V_1\otimes\mathbb{C}_k),\quad V_6\oplus V_{8,i}\oplus V_{5,j}\oplus\mathbb{C}_k,
$$
and
$$
  V_6\oplus(V_{8,i}\otimes_{\mathbb{C}}\mathbb{C}_k)\oplus V_{5,j}\oplus(V_1\otimes\mathbb{C}_k),\quad V_6\oplus(V_{8,i}\otimes_{\mathbb{C}}\mathbb{C}_k)\oplus V_{5,j}\oplus\mathbb{C}_k.
$$
Consider $C_2\trianglelefteq\TL(2,5)$ and the action of the quotient $S_5\times C_k$ on the $\mathbb{Z}_2$-homology sphere $\Sigma'=\Sigma^{C_2}$. In our case, the tangent space at the fixed point $x$ of this action is isomorphic to an $\mathbb{R}(S_5\times C_k)$-representation of one of the forms:
$$
V_6\oplus\mathbb{C}_k,\quad
V_6\oplus(V_1\otimes\mathbb{C}_k),\quad
V_6\oplus V_{5,i}\oplus\mathbb{C}_k,\quad
V_6\oplus V_{5,i}\oplus(V_1\otimes \mathbb{C}_k).
$$
We now apply \Cref{corolarry:one_dims_struct} with
$$
G=S_5,\qquad H=A_5,\qquad G/H\cong C_2.
$$
Here $V_1$ is the non-trivial one-dimensional representation of $S_5/A_5\cong C_2$, while $V_6$ and $V_{5,i}$ are regarded as $S_5$-representations with trivial $C_k$-action. Thus the above representations have the form required in \Cref{corolarry:one_dims_struct}. Consequently, the fixed point set of the induced $S_5\times C_k$-action on $\Sigma'$ is not a singleton. Therefore
$$
\Sigma^{\TL(2,5)\times C_k}=(\Sigma')^{S_5\times C_k}\neq\{\pt\},
$$
which gives the desired contradiction.
\end{proof}

\begin{lemma}\label{lemma:structural_s5}
Let $W$ be one of the following representations of $S_{5}\times C_k$:
\begin{enumerate}
    \item $   V_{4,2}^{\oplus k_{1}}   \oplus V_{5,1}^{\oplus k_{2}}   \oplus V_{5,2}^{\oplus k_{3}} \oplus\mathbb{C}_k^{\oplus k_{4}} \oplus (V_{1}\otimes \mathbb{C}_k)^{\oplus k_{5}} \oplus (V_{4,2}\otimes \mathbb{C}_k)^{\oplus k_{6}}\oplus (V_{5,1}\otimes \mathbb{C}_k)^{\oplus k_{7}}  $,
 
    \item  $V_{4,1}^{\oplus k_{1}}   \oplus V_{4,2}^{\oplus k_{2}} \oplus V_{5,1}^{\oplus k_{3}}\oplus V_{5,2}^{\oplus k_{4}}  \oplus\mathbb{C}_k^{\oplus k_{5}} \oplus (V_{1}\otimes \mathbb{C}_k)^{\oplus k_{6}}\oplus (V_{4,1}\otimes \mathbb{C}_k)^{\oplus k_{7}}\oplus (V_{5,2}\otimes \mathbb{C}_k)^{\oplus k_{8}}     $.

\end{enumerate}
Then $W$ cannot occur as the tangential
representation at the fixed point of a one fixed point action
on a homology sphere where $k_{i}$
are arbitrary non-negative integers.

\end{lemma}

\begin{proof}
 In the first two cases, we shall proceed similarly to the proofs of \cref{lemma:general_form_sl25,lemma:structural_tl25}.
    \begin{enumerate}
        \item For subgroups $ C_{4}\leq (C_{5}\rtimes C_{4}) \times C_{k} \cap D_{8} $ 
        we have the following computations of the dimensions of the fixed point subspaces:
        $$
        \mathrm{dim}W^{ C_{4}}=  k_{1}   + k_{2} + 2k_{3} +2k_{4}+ 2k_{6}+2k_{7}, $$
         $$ \mathrm{dim}W^{D_{8}}=  k_{1}  + k_{2} + k_{3} +2k_{4}+ 2k_{6}+2k_{7}, $$
        $$ \mathrm{dim}W^{( C_{5}\rtimes C_{4}) \times C_{k}  } =  k_{3},  $$
       and the equality 
         $\mathrm{dim}W^{ C_{4}} = \mathrm{dim}W^{D_{8}} +\mathrm{dim}W^{(C_{5}\rtimes C_{4}) \times C_{k} }$
         holds. Since $D_{8}$ and $(C_{5}\rtimes C_{4}) \times C_{k}$ 
         generate $S_{5}\times C_k$, the result follows.
        \item We use the same argument as in case 1
        for subgroups $(C_2^2)_{15}, D_{8}$ and $D_{12}\times C_{k}$: 
        $$
            \mathrm{dim}W^{ (C_2^2)_{15}}=  2k_{2} + 2k_{3} +k_{4}+ 2k_{5}+2k_{8},$$
         $$ \mathrm{dim}W^{D_{12}\times C_k}= k_{2}+k_{3},$$
        $$ \mathrm{dim}W^{D_{8}} =  k_{2} +k_{3} + k_{4}+ 2k_{5}+2k_{8} .$$
     
    \end{enumerate}

\end{proof}

\begin{theorem}\label{theorem:3main2}\emph{(cf. \Cref{theorem:3main})}
$S_{5}\times C_k$ does not admit a one fixed point action on a homology sphere $\Sigma$ of dimension $n$ when $n\in\{6,\ldots,13\}\cup\{15\}$.  
 
\end{theorem}

\begin{proof}
    Let $V$ be the tangential representation at the single fixed point $x$. By \Cref{lemma:structural_s5} and \Cref{lemma:one_dim} it follows that $V$ must contain at least one representation from each of the sets below as a direct summand: 
    $$
    \mathcal{S}_1=\{V_{4,1},V_{4,1}\otimes\mathbb{C}_k,V_{5,2}\otimes\mathbb{C}_k,V_6,V_6\otimes\mathbb{C}_k\},
    $$
    $$
    \mathcal{S}_2=\{V_{4,2}\otimes\mathbb{C}_k,V_{5,1}\otimes\mathbb{C}_k,V_6,V_6\otimes\mathbb{C}_k\}.
    $$
    
    \noindent Moreover, if $V$ is odd-dimensional, then by the dimension constraint and \Cref{lemma:one_dim}, $V$ must also contain $V_{5,1}$ or $V_{5,2}$. Moreover, by \Cref{corolarry:one_dims_struct}, we can exclude the case when $V$ does not contain a summand of the form $W\otimes\mathbb{C}_k$, where $W$ is not $V_1$. Therefore, combining all of the above information with \Cref{table:dims_s5}, we are left only with the dimension $12$ and the modules:
    $$
    V\cong V_{6}\otimes \mathbb{C}_k\quad\text{and}\quad V\cong V_{4,1}\oplus(V_{4,2}\otimes\mathbb{C}_k).
    $$

Let $V=V_{4,1}\oplus(V_{4,2}\otimes\mathbb{C}_k)$. Take
$$
H=C_5,\quad H_1=C_5\times C_k,\quad H_2=S_5.
$$
Then $H\leq H_1\cap H_2$, the subgroups $H_1$ and $H_2$ generate $S_5\times C_k$, and
$$
\dim V^H=0.
$$
The subgroup $H=C_5$ is cyclic of prime order, hence it is an extension of a cyclic group by a $5$-group, and all its elements have prime-power order. Moreover, both $H_1=C_5\times C_k$ and $H_2=S_5$ cannot act with exactly one fixed point in the required tangential representation type. For $H_1$, this follows from the non-Oliver sequence
$$
C_5\trianglelefteq C_5\times C_k\trianglelefteq C_5\times C_k.
$$
For $H_2=S_5$, this follows from \cite[Theorem 4.1]{MorimotoTamura}, since the restricted tangential representation is isomorphic to
$$
\mathrm{Res}_{S_5}^{S_5\times C_k}(V)\cong V_{4,1}\oplus V_{4,2}^{\oplus 2}.
$$
Thus the hypotheses of \Cref{lemma:equal_dims} are satisfied, and \Cref{lemma:equal_dims} gives a contradiction.


In the case $V= V_{6}\otimes \mathbb{C}_k$, we observe that there is a cyclic subgroup of prime power order $P\leq C_k$ acting non-trivially on $V$. Hence, $\Sigma^{S_5\times C_k}=(\Sigma^P)^{(S_5\times C_k)/P}=(S^0)^{(S_5\times C_k)/P}=S^0$.

    \end{proof}

We are also able to obtain the following exclusion results concerning dimension $17$.
\begin{lemma}\label{lemma:37}
    Let $S_5 \times C_{k}$ act on a $17$-dimensional homology sphere $\Sigma$ with a fixed point $x$ and tangential $S_5 \times C_{k}$-representation $T_{x}\Sigma$ isomorphic to $V_{5,i}\oplus V_{4,1}\oplus(V_{4,2}\otimes\mathbb{C}_k)$ for $i=1,2$. Then $\Sigma^{S_5 \times C_{k}}$ contains an even number of points.

\end{lemma}

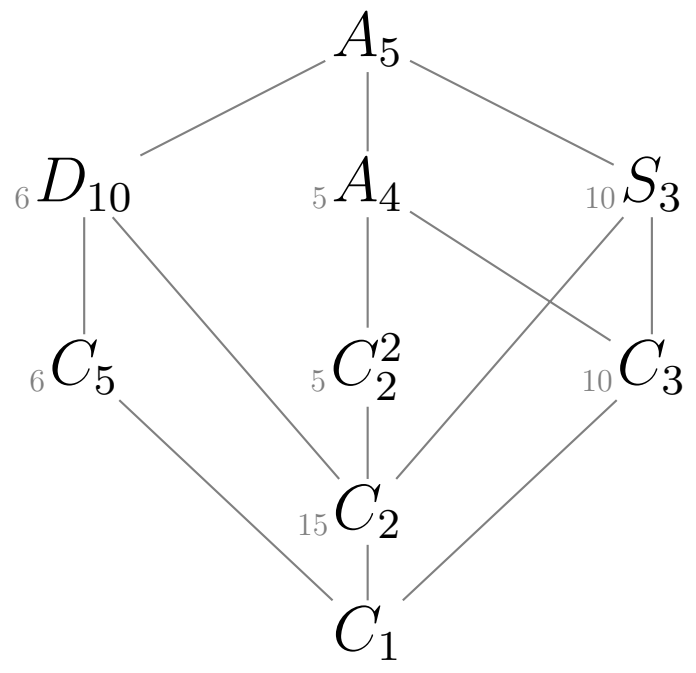
\begin{figure}[H]
  \centering
   \begin{tikzpicture}[scale=1.0,sgplattice]
\node at (2,0) (1) {\gn{C1}{C_1}};
  \node at (2,0.803) (2) {\gn{C2}{C_2}};
  \node at (3.88,1.76) (3) {\gn{C3}{C_3}};
  \node at (0.125,1.76) (4) {\gn{C5}{C_5}};
  \node at (2,1.76) (5) {\gn{C2^2}{C_2^2}};
  \node at (3.88,2.97) (6) {\gn{S3}{S_3}};
  \node at (0.125,2.97) (7) {\gn{D5}{D_{10}}};
  \node at (2,2.97) (8) {\gn{A4}{A_4}};
  \node at (2,3.93) (9) {\gn{A5}{A_5}};
  \draw[lin] (1)--(2) (1)--(3) (1)--(4) (2)--(5) (2)--(6) (3)--(6) (2)--(7)
     (4)--(7) (3)--(8) (5)--(8) (6)--(9) (7)--(9) (8)--(9);
  \node[cnj=2] {15};
  \node[cnj=3] {10};
  \node[cnj=4] {6};
  \node[cnj=5] {5};
  \node[cnj=6] {10};
  \node[cnj=7] {6};
  \node[cnj=8] {5};
\end{tikzpicture}
 \caption{Conjugacy classes subgroups of $A_5$.}
  \label{fig:A5}
\end{figure}

\begin{proof}
  First we show that the action restricted to $S_5$ has an even number of fixed points, and from this we deduce that the same holds for $S_5\times C_k$. Every element of $A_5$ has order $1$, $2$, $3$, or $5$. Hence every element of $A_5$ has prime-power order, and every element of $2$-power order has order at most $2$, in particular at most $4$. Thus the $8$-condition required in \Cref{lemma:8cond} is satisfied. Therefore, by \Cref{lemma:8cond}, $\Sigma^{A_5}$ is a discrete set. Consider the restriction of the representation $\mathrm{Res}^{S_5\times C_k}_{A_5}(V) \cong V_{4}^{\oplus 3} \oplus V_{5}$ (see \Cref{table:fpsdim_A5}). By the $\mathrm{mod}$ $2$ intersection number of $\Sigma^{C^{2}_{2}}$ and $\Sigma^{S_3}$ in $\Sigma^{C_2}$, $\Sigma^{A_5}$ contains an even number of points (the appropriate subgroups can be consulted in \Cref{fig:A5}), and the same holds for $\Sigma^{S_5}$, since $S_{5}/A_{5} \cong C_{2}$.

  To show the same conclusion for $S_{5} \times C_{k}$ we use the subgroups $C_{5} < C_{5} \times C_{k}$ and $S_{5}$ together with the fact that $\Sigma^{C_{5}} = \Sigma^{C_{5} \times C_{k}}$ and $\Sigma^{S_{5}}$ contains an even number of points by the first part of the proof -- it suffices to observe $\Sigma^{S_5}=\Sigma^{S_5\times C_k}$.
\end{proof}

\begin{lemma}\label{lemma:3_8}\emph{(cf. \Cref{theorem:3main})}
If $k=3l$ then the action of $S_{5}\times C_{k}$ on a $17$-dimensional homology sphere with the tangential representation at one of the fixed points isomorphic to $(V_{6}\otimes \mathbb{C}_{k}) \oplus V_{5,i}$ has exactly two fixed points.
    
\end{lemma}

\begin{proof}
    Suppose $k=3l$ and let $S_5\times C_k$ act on a $17$-dimensional homology sphere $\Sigma$. Let $\mathcal{C}_3\leq C_k$ denote the subgroup of order $3$ and $V=(V_{6}\otimes \mathbb{C}_{k}) \oplus V_{5,i}$ be the tangential representation at one of the fixed points, denoted by $x$.

    First we show that $\Sigma^{A_4}=S^0$. By \Cref{table:dims_s5}, for the representation
    $V=(V_{6}\otimes \mathbb{C}_{k}) \oplus V_{5,i}$ we have
    $$
    \dim V^{(C_2^2)_5}=2
    \quad\text{and}\quad
    \dim V^{A_4}=0.
    $$
    By Smith theory, $\Sigma^{(C_2^2)_5}$ is a $2$-dimensional $\mathbb{Z}_2$-homology sphere,
    hence, by the classification of surfaces, it is homeomorphic to $S^2$. Since
    $(C_2^2)_5\trianglelefteq A_4$ and $A_4/(C_2^2)_5\cong C_3$, we get
    $$
    \Sigma^{A_4}
    =
    \left(\Sigma^{(C_2^2)_5}\right)^{A_4/(C_2^2)_5}.
    $$
    The fixed point set on the right-hand side is a $\mathbb{Z}_3$-homology sphere by Smith theory, and it is
    discrete because $\dim V^{A_4}=0$. Therefore
    $$
    \Sigma^{A_4}\cong S^0.
    $$

    To show that $\Sigma^{S_5\times C_k}$ consists of two points, it is enough to prove that
    $\Sigma^{A_5\times \mathcal{C}_3}=S^0$ and then use
    $$
    \Sigma^{S_5\times C_k}
    =
    \left(\Sigma^{A_5\times \mathcal{C}_3}\right)^{C_2\times C_l}.
    $$
    We first note that
    $
    \Sigma^{A_4\times \mathcal{C}_3}=\Sigma^{A_4}=S^0.
    $
    Next we justify the equality
    $$
    \Sigma^{C_3\times \mathcal{C}_3}
    =
    \Sigma^{S_3\times \mathcal{C}_3},
    $$
    where $C_3\leq S_3\leq A_5$. Indeed, from the fixed point dimension table and by the non-triviality of the action of $\mathcal{C}_3$ (which follows from the effectiveness of the action of $S_5\times C_k$), we have
    $$
    \dim V^{C_3\times \mathcal{C}_3}
    =
    \dim V^{S_3\times \mathcal{C}_3}=1.
    $$
    Hence the connected components of $\Sigma^{C_3\times \mathcal{C}_3}$ and
    $\Sigma^{S_3\times \mathcal{C}_3}$ containing $x$ have the same dimension. By the Principal Orbit
    Theorem, the latter fixed point set is open and dense in the former one. Since it is
    also a closed submanifold of $\Sigma^{C_3\times \mathcal{C}_3}$, which is connected as a positive-dimensional $\mathbb{Z}_3$-homology sphere, the two fixed point
    submanifolds are equal:
    $$
    \Sigma^{C_3\times \mathcal{C}_3}
    =
    \Sigma^{S_3\times \mathcal{C}_3}.
    $$

    \noindent Finally, since $A_5$ is generated by $A_4$ and $S_3$, we obtain
    $$
    \Sigma^{A_5\times \mathcal{C}_3}
    =
    \Sigma^{A_4\times \mathcal{C}_3}
    \cap
    \Sigma^{S_3\times \mathcal{C}_3}
    =
    \Sigma^{A_4\times \mathcal{C}_3}
    \cap
    \Sigma^{C_3\times \mathcal{C}_3}
    =
    \Sigma^{A_4\times \mathcal{C}_3}
    =
    S^0.
    $$
    Consequently,
    $$
    \Sigma^{S_5\times C_k}
    =
    \left(\Sigma^{A_5\times \mathcal{C}_3}\right)^{C_2\times C_l}
    =
    S^0.
    $$
\end{proof}

\subsection{The case $k=2$}\label{subsection:3_2}
Let us deal first with the case $G\times C_{2}$ for $G\in\{\SL(2,5),\TL(2,5),S_5\}$. As we shall see, this case is much simpler than the previous case with an odd-order cyclic group as a direct summand. The reason for this is that we can make better use of $2$-group exclusion techniques. The following lemma is crucial.

\begin{lemma}\label{lemma:c2_general_structure}
If   $\mathrm{SL}(2,5)\times C_{2}$ (resp. $\mathrm{TL}(2,5)\times C_{2}$,  $S_{5}\times C_{2}$) acts on a homology sphere $\Sigma$ with exactly one fixed point $x$, then the tangential representation at $x$ contains a $3$-dimensional (resp. $6$-dimensional) irreducible summand $V_3$ (resp. $V_{6}$).  
 
\end{lemma}

\begin{proof}

We distinguish triples of subgroups for each of the groups in question. In each case the group in the first column is a $2$-group, and the groups appearing in the second and the third column contain that $2$-group. Moreover, they generate the whole group in question and none of them is an Oliver group.

We consider only the non-trivial irreducible representations, except $V_3$ in the case of $\SL(2,5)\times C_2$ and except $V_6$ in the cases of $\TL(2,5)\times C_2$ and $S_5\times C_2$, for which the fixed point dimensions are not all zero for the corresponding subgroup triple. If all these fixed point dimensions are zero, then we apply \Cref{lemma:equal_dims} directly.

 \begin{table}[H]
\begin{center}
{\fontsize{10pt}{12pt}\selectfont

\begin{tabular}{|c||c|c|c|} \hline
$H$ & 

 $C_{4} \times C_{2}$&$Q_{8} \times C_{2}$ & 
$(C_{3}\rtimes C_{4}) \times C_{2}$ \\ 
\hline
\hline

$V_{4,1}$ & 
 $2$ &$1$ & $1$ \\
\hline

$V_{5}$ & 
 $3$ &$2$ & $1$  \\

\hline

\end{tabular}
}
\vskip 0.2cm
\caption{Fragment of \Cref{table:dims_sl25} of our interest.}\label{table:c2_table1}
\end{center}
\end{table}

    \begin{table}[H]
\begin{center}
    {\fontsize{10pt}{12pt}\selectfont 

    \begin{tabular}{|c||c|c|c|} \hline
$H$ & 

 $(Q_{8})_{15}\times C_{2}$&$Q_{16} \times C_{2}$ & 
$(C_3 \rtimes Q_{8}) \times C_{2}$ \\ 
\hline
\hline

$V_{4,2}$ & 
 $2$ &$1$ & $1$ \\
\hline

$V_{5,1}$ & 
 $2$ &$1$ & $1$  \\
 \hline

$V_{5,2}$ & 
 $1$ &$1$ & $0$ \\
\hline

\end{tabular}
}
\vskip 0.2cm
\caption{Fragment of \Cref{table:dims_tl25} of our interest.}\label{table:c2_table2}
\end{center}
\end{table}

\begin{table}[H]
\begin{center}
{\fontsize{10pt}{12pt}\selectfont

\begin{tabular}{|c||c|c|c|} \hline
$H$ & 

 $(C^{2}_{2})_{15}\times C_{2}$&$D_{8} \times C_{2}$ & 
$D_{12} \times C_{2}$ \\ 
\hline
\hline

$V_{4,2}$ & 
 $2$ &$1$ & $1$ \\
\hline

$V_{5,1}$ & 
 $2$ &$1$ & $1$ \\
\hline 
$V_{5,2}$ & 
 $1$ &$1$ & $0$  \\
 \hline

\end{tabular}
}
\vskip 0.2cm
\caption{Fragment of \Cref{table:dims_s5} of our interest.}\label{table:fpsdim_V3}
\end{center}
\end{table}

\noindent We now explain how the data in the tables are used. Let $W$ be a direct sum of irreducible summands appearing in one of the corresponding tables, and let
$$
H\leq H_1\cap H_2
$$
be the subgroup triple displayed in that table. 

If, for one of $i=1,2$, we have
$
\dim W^H=\dim W^{H_i}>0,
$
then we apply \Cref{theorem:equaldims_Kamil}. 
In the case
$$
\dim W^H=\dim W^{H_1}+\dim W^{H_2}
$$
with
$
\dim W^{H_1}>0
\quad\text{and}\quad
\dim W^{H_2}>0,
$
then we apply \Cref{theorem:pdis}.

The two alternatives above, together with the case $\dim W^H=0$ discussed before, cover all possible direct sums of the representations displayed in the tables. Hence, if a one fixed point action exists, its tangential representation must contain the exceptional summand not covered by the tables, namely $V_3$ in the case of $\SL(2,5)\times C_2$, and $V_6$ in the cases of $\TL(2,5)\times C_2$ and $S_5\times C_2$.

The proof is concluded, we explain, however, why summands of the form $U\otimes\mathbb{R}_{\pm}$ do not need to be treated separately. All subgroup triples used above contain the second factor $C_2$. Since $\mathbb{R}_{\pm}$ is non-trivial on this factor, we have
$$
(U\otimes\mathbb{R}_{\pm})^K=0
$$
for every subgroup $K$ appearing in the tables. Thus such summands do not change any of the fixed point dimension computations above. 
  
\end{proof}

\begin{theorem}\emph{(cf. \Cref{theorem:1main})}\label{theorem:1main2only}
 $\mathrm{SL}(2,5)\times C_{2}$ does not admit a one fixed point action on a homology sphere $\Sigma$ of dimension $n$ when $n\in\{7,8,9\}$.
  
\end{theorem}
\begin{proof}
By \Cref{lemma:c2_general_structure} we know that the representation $V$ appearing at the fixed point has to contain $V_3$. On the other hand, by \Cref{lemma:one_dim}, $V$ must not contain a one-dimensional irreducible summand. Note that all the irreducible representations of $\SL(2,5)$ except $V_{4,2}$, $V_8$ and $V_{12}$ have kernel equal to $C_2$. Thus, the restriction of $V$ to $\SL(2,5)$ has to contain one of these representations as a direct summand. On the other hand, $C_{2}$ (more precisely, $\{1\}\times C_{2}$) acts effectively on $V$, so $V$ must also contain a direct summand of the form $V\otimes \mathbb{R}_{\pm}$ for some representation $V$ of $\mathrm{SL}(2,5)$. These constraints leave potentially the case $V=V_3\oplus(V_{4,2}\otimes\mathbb{R}_{\pm})$. The diagonal subgroup (isomorphic to the cyclic group of order two) generated by $(a,b)$, where $a$ is the generator of $C_2$ in $\SL(2,5)$ and $b$ is the generator of $\{1\}\times C_2$, acts trivially on $V$, however. Hence, no admissible faithful representation in dimensions $n\in\{7,8,9\}$ can occur. This concludes the proof.


\end{proof}

\begin{theorem}\emph{(cf. \Cref{theorem:2main})}\label{theorem:2main2only}
    Let $G=\TL(2,5)\times C_{2}$. Then there does not exist a one fixed point action of $G$ on $S^n$ whenever $n\in\{6,\ldots,17\}\cup\{21\}$.
\end{theorem}
\begin{proof}
    We proceed similarly to the proof of the previous theorem. All the irreducible $\TL(2,5)$-representations except $V_{8,i}$ and $V_{12}$ have kernel equal to $C_{2}$. By the effectiveness assumption, the restriction of the representation $V$ at the tangent space at the fixed point to $\TL(2,5)$ has to contain one of these representations as a direct summand. On the other hand, by \Cref{lemma:c2_general_structure}, $V$ has to contain $V_6$ as a direct summand as well.

    Using the dimensions of the representations appearing in \Cref{table:dims_tl25}, together with the effectiveness of the restriction to the second factor $\{1\}\times C_2$, the only possible remaining representation is
$$
V\cong (V_{8,i}\otimes\mathbb{R}_{\pm})\oplus V_6.
$$
Indeed, the summand $V_6$ is forced by \Cref{lemma:c2_general_structure}, while the effectiveness of the $\{1\}\times C_2$-action forces the presence of a summand on which this second factor acts non-trivially. In the dimensions under consideration, the only such remaining possibility is $V_{8,i}\otimes\mathbb{R}_{\pm}$.

We note, however, that $V$ is not faithful. Namely, $V_{8,i}\otimes\mathbb{R}_{\pm}$, in contrast to the representation $V_{8,i}$ itself, is not a faithful $\TL(2,5)\times C_2$-representation, and its kernel is contained in the kernel of $V_6$. Hence the whole representation $V$ is not faithful. This contradicts the effectiveness of the original action and concludes the proof.
    

\end{proof}

\begin{theorem}\emph{(cf. \Cref{theorem:3main})}\label{theorem:3main2only}
    Let $G=S_5\times C_{2}$. Then there does not exist a one fixed point action of $G$ on $S^n$ whenever $n\in\{6,\ldots,9\}\cup\{13\}$.
\end{theorem}
\begin{proof}

Let $V$ be the representation at the single fixed point. By \Cref{lemma:c2_general_structure}, $V$ must contain $V_6$ as a direct summand. Moreover, since the fixed point is isolated, $V$ cannot contain the trivial summand, and by \Cref{lemma:one_dim} it cannot contain a non-trivial one-dimensional summand. Finally, the effectiveness of the action of the second factor, i.e. $\{1\}\times C_2$, implies that $V$ must contain a summand of the form $U\otimes\mathbb{R}_{\pm}$.

If $n\leq 9$, then after the forced summand $V_6$ there remain at most three dimensions. Since one-dimensional summands are excluded and the smallest non-one-dimensional irreducible $S_5$-representation has dimension $4$, there is no room for a summand on which the second factor $C_2$ acts non-trivially. Hence no faithful tangential representation can occur.

Now let $n=13$. Again $V$ contains $V_6$, and it must also contain a summand $U\otimes\mathbb{R}_{\pm}$ with $U$ not one-dimensional. By \Cref{table:dims_s5}, the possible dimensions of such a non-one-dimensional irreducible $S_5$-representation are at least $4$. If $\dim U=4$, then three dimensions remain; if $\dim U=5$, then two dimensions remain; and if $\dim U=6$, then one dimension remains. None of these remaining dimensions can be filled by admissible irreducible summands, since the trivial and non-trivial one-dimensional summands are excluded and there are no irreducible $S_5$-representations of dimensions $2$ or $3$. Therefore no faithful tangential representation exists in dimension $13$ either.
\end{proof}

\subsection{The case of even $k>2$}
\begin{lemma}\label{lemma:general_form_sl25_even}
The following representation $W$ of $\SL(2,5)\times C_k$:

     $$    V_{4,1}^{\oplus k_{1}}   \oplus V_{4,2}^{\oplus k_{2}} \oplus V_{5}^{\oplus k_{3}}   \oplus V_{8}^{\oplus k_{4}}   \oplus V_{12}^{\oplus k_{5}} \oplus \mathbb{C}_k^{\oplus k_{6}}\oplus (V_{4,2}\otimes_{\mathbb{C}} \mathbb{C}_k)^{\oplus k_{7}}\oplus  (W'\otimes\mathbb{R}_{\pm}),$$
     where $W'$ is some real $\SL(2,5)$-representation, cannot occur as the tangential
representation at the fixed point of a one fixed point action
on a homology sphere where $k_{i}$
are arbitrary non-negative integers.

\end{lemma}

\begin{proof}
The proof is analogous to the proof of \Cref{lemma:general_form_sl25}, but we spell out the treatment of the summands involving $C_k$. Write
$$
C_k\cong C_{2^l}\times C_m,
$$
where $m$ is odd. We consider the subgroups
$$
H=C_4\times C_{2^l},\quad
H_1=Q_8\times C_{2^l},\quad
H_2=(C_3\rtimes C_4)\times C_k.
$$
Then $H\leq H_1\cap H_2$ and
$$
\langle H_1,H_2\rangle=\SL(2,5)\times C_k.
$$

\noindent The summands of the form $W'\otimes\mathbb{R}_{\pm}$ do not contribute to the fixed point dimensions for these subgroups, since each of the subgroups $H$, $H_1$, and $H_2$ contains a subgroup on which $\mathbb{R}_{\pm}$ acts non-trivially.

The only additional point, compared with \Cref{lemma:general_form_sl25}, is that a summand $\mathbb{C}_k$ need not have zero fixed points after restriction to $C_{2^l}$. We therefore put
$$
d=\dim\left(\mathbb{C}_k^{\oplus k_6}\right)^{C_{2^l}}.
$$
Applying \Cref{table:dims_sl25} and \Cref{corollary:fixed_dims2}, we obtain
$$
\mathrm{dim}W^{H}=2k_1+3k_3+d,
$$
$$
\mathrm{dim}W^{H_1}=k_1+2k_3+d,
$$
and
$$
\mathrm{dim}W^{H_2}=k_1+k_3.
$$
Thus
$$
\mathrm{dim}W^{H}=\mathrm{dim}W^{H_1}+\mathrm{dim}W^{H_2}.
$$

\noindent We now distinguish three cases.

First assume that $k_1=k_3=d=0$. Then $\mathrm{dim}W^{H}=0$. We apply \Cref{lemma:equal_dims} to the subgroups
$$
H=C_4\times C_{2^l},\quad
H_1=Q_8\times C_{2^l},\quad
H_2=(C_3\rtimes C_4)\times C_k.
$$
The subgroup $H$ is a $2$-group, hence it is an extension of a cyclic group by a $2$-group, and all its elements have prime-power order. Moreover, $H_1$ and $H_2$ generate $\SL(2,5)\times C_k$, and $\mathrm{dim}W^{H}=0$. Finally, both $H_1$ and $H_2$ cannot act with exactly one fixed point in the required tangential representation type. For $H_2$, this follows from the non-Oliver sequence
$$
C_3\trianglelefteq C_3\times C_k\trianglelefteq (C_3\rtimes C_4)\times C_k.
$$
Thus \Cref{lemma:equal_dims} gives a contradiction.

Now assume that $k_1=k_3=0$ and $d>0$. Then
$$
\mathrm{dim}W^H=d=\mathrm{dim}W^{H_1}.
$$
We apply \Cref{theorem:equaldims_Kamil} to the subgroups
$$
H=C_4\times C_{2^l},\quad
H_1=Q_8\times C_{2^l},\quad
H_2=(C_3\rtimes C_4)\times C_k.
$$
The subgroup inclusions and the generation condition are as above, and the equality $\mathrm{dim}W^{H_1}=\mathrm{dim}W^H$ holds. Since $H_2=(C_3\rtimes C_4)\times C_k$ cannot act with exactly one fixed point in the required tangential representation type, \Cref{theorem:equaldims_Kamil} gives a contradiction.

Finally assume that $k_1+k_3>0$. Then
$$
\mathrm{dim}W^{H_1}>0
$$
and
$$
\mathrm{dim}W^{H_2}>0.
$$
We apply \Cref{theorem:pdis} to the subgroups
$$
H=C_4\times C_{2^l},\quad
H_1=Q_8\times C_{2^l},\quad
H_2=(C_3\rtimes C_4)\times C_k.
$$
The subgroup inclusions and the generation condition are as above. The fixed point manifolds appearing in \Cref{theorem:pdis} are orientable in this case, and the equality
$$
\mathrm{dim}W^{H}=\mathrm{dim}W^{H_1}+\mathrm{dim}W^{H_2}
$$
gives the required complementary-dimension condition. Therefore \Cref{theorem:pdis} gives a contradiction.

In all cases we obtain a contradiction, and hence $W$ cannot occur as the tangential representation at the fixed point of a one fixed point action on a homology sphere.

\end{proof}

\begin{theorem}\label{theorem:1main2_even}\emph{(cf. \Cref{theorem:1main})}
    $\SL(2,5)\times C_k$ does not admit a one fixed point action on $S^n$ if $n\in\{7,8,9\}$.
\end{theorem}
\begin{proof}


Assume, for a contradiction, that such an action exists. Let $\Sigma$ be a homology
sphere with a one-fixed-point action of $\SL(2,5)\times C_k$, and let $x$ be its unique
fixed point. Denote by
$
        V=T_x\Sigma
$ the tangential representation at the fixed point.

In listing the possible tangential representations below, we use simultaneously
\Cref{table:dims_sl25}, \Cref{corollary:fixed_dims2}, the condition
$
        V^{\SL(2,5)\times C_k}=0,
$
the exclusion of one-dimensional summands by \Cref{lemma:one_dim}, and the faithfulness of the
$\SL(2,5)\times C_k$-representation $V$. The only candidates which are not covered directly by \Cref{lemma:general_form_sl25_even} remain in dimension $9$ and are


$$
        V\cong V_3\oplus V_{4,2}\oplus \mathbb C_k
$$
and
$$
        V\cong V_3\oplus (V_{4,2}\otimes_{\mathbb C}\mathbb C_k)\oplus \mathbb C_k.
$$
We now exclude these two cases.

Let
$$
        Z=Z(\SL(2,5))\cong C_2.
$$
Then $Z$ is normal in $\SL(2,5)$ and
$$
        \SL(2,5)/Z\cong A_5.
$$
Since $Z$ is normal in $\SL(2,5)\times C_k$, the connected component $\Sigma^Z_x$ is
invariant under the action of $\SL(2,5)\times C_k$, and the induced action on
$\Sigma^Z_x$ factors through
$$
        (SL(2,5)\times C_k)/Z\cong A_5\times C_k.
$$

\noindent For both representations listed above, we have
$$
        V^Z\cong V_3\oplus \mathbb C_k.
$$
Indeed, the summands $V_{4,2}$ and
$V_{4,2}\otimes_{\mathbb C}\mathbb C_k$ have no $Z$-fixed vectors, whereas $Z$ acts
trivially on $V_3$ and on the pure $\mathbb C_k$-summand. Hence
$$
        T_x\Sigma^Z_x\cong V^Z\cong V_3\oplus \mathbb C_k
$$
as an $A_5\times C_k$-representation.

We apply \Cref{corolarry:one_dims_struct} to the induced $A_5\times C_k$-action on $\Sigma^Z_x$, using
the pair of subgroups $A_5$ and $C_k$. For the representation $V^Z$ we have from \Cref{table:fpsdim_A5}
$$
        \dim (V^Z)^{A_5}=2
        \qquad\text{and}\qquad
        \dim (V^Z)^{C_k}=3.
$$
Thus
$$
        \dim (V^Z)^{A_5}+\dim (V^Z)^{C_k}
        =
        5
        =
        \dim V^Z.
$$
Moreover,
$$
        (V^Z)^{A_5\times C_k}=0.
$$
Therefore the hypotheses of \Cref{corolarry:one_dims_struct} are satisfied, and we conclude that
$$
        (\Sigma^Z_x)^{A_5\times C_k}
$$
cannot be a singleton. But
$$
        (\Sigma^Z_x)^{A_5\times C_k}
        =
        \Sigma^{\SL(2,5)\times C_k}.
$$
This contradicts the assumption that the original action has exactly one fixed point.
Therefore the two remaining $9$-dimensional candidates are excluded.
\end{proof}

\begin{lemma}\label{lemma:structural_tl25_even}
Let $W$ be one of the following representations of $\TL(2,5)\times C_k$:
\begin{enumerate}
    \item $ V_{4,1}^{\oplus k_{1}}  \oplus V_{4,2}^{\oplus k_{2}} \oplus V_{5,1}^{\oplus k_{3}}   \oplus V_{5,2}^{\oplus k_{4}}   \oplus V_{8,1}^{\oplus k_{5}} \oplus V_{8,2}^{\oplus k_{6}}  \oplus V_{12}^{\oplus k_{7}}\oplus\mathbb{C}_k^{\oplus k_{8}} \oplus (V_{1}\otimes \mathbb{C}_k)^{\oplus k_{9}}\oplus (V_{4,1}\otimes \mathbb{C}_k)^{\oplus k_{10}} \oplus (V_{5,2}\otimes \mathbb{C}_k)^{\oplus k_{11}} \oplus (V_{8,1}\otimes_{\mathbb{C}}\mathbb{C}_k)^{k_{12}}\oplus(V_{8,2}\otimes_{\mathbb{C}}\mathbb{C}_k)^{k_{13}}\oplus (V_{12}\otimes_{\mathbb{C}}\mathbb{C}_k)^{k_{14}}\oplus (W'\otimes\mathbb{R}_{\pm})  $
    \item  $V_{4,2}^{\oplus k_1}
\oplus V_{5,1}^{\oplus k_2}
\oplus V_{5,2}^{\oplus k_3}
\oplus V_{8,1}^{\oplus k_4}
\oplus V_{8,2}^{\oplus k_5}
\oplus V_{12}^{\oplus k_6}
\oplus \mathbb C_k^{\oplus k_7} 
\oplus (V_1\otimes \mathbb C_k)^{\oplus k_8}
\oplus (V_{4,2}\otimes \mathbb C_k)^{\oplus k_9}
\oplus (V_{5,1}\otimes \mathbb C_k)^{\oplus k_{10}} 
\oplus (V_{8,1}\otimes_{\mathbb C}\mathbb C_k)^{\oplus k_{11}}
\oplus (V_{8,2}\otimes_{\mathbb C}\mathbb C_k)^{\oplus k_{12}}
\oplus (V_{12}\otimes_{\mathbb C}\mathbb C_k)^{\oplus k_{13}}
\oplus (W'\otimes \mathbb R_{\pm})$
\end{enumerate}
where $W'$ is a real $\TL(2,5)$-representation. Then $W$ cannot occur as the tangential
representation at the fixed point of a one fixed point action
on a homology sphere where $k_{i}$
are arbitrary non-negative integers.

\end{lemma}



\begin{proof}
Assume, for a contradiction, that there exists a homology sphere $\Sigma$ with a
one-fixed-point action of $\TL(2,5)\times C_k$ and with tangential representation at the
fixed point isomorphic to $W$. Write
$$
        C_k\cong C_{2^\ell}\times C_m,
$$
where $m$ is odd. Since $k$ is even, we have $\ell\geq 1$.

We first consider case $(1)$. We use the subgroups
$$
        H={Q_8}_{(15)}\times C_{2^\ell},\qquad
        H_1=Q_{16}\times C_{2^\ell},\qquad
        H_2=(C_3\rtimes Q_8)\times C_k.
$$
Then
$$
        H\leq H_1\cap H_2
$$
and
$$
        \langle H_1,H_2\rangle=\TL(2,5)\times C_k.
$$
Using \Cref{table:dims_tl25} and \Cref{corollary:fixed_dims2}, we have
$$
        \dim W^{H_2}=k_2+k_3.
$$
Moreover,
$$
        \dim W^H=\dim W^{H_1}+k_2+k_3.
$$
Indeed, the summands $V_{4,2}^{\oplus k_2}$ and $V_{5,1}^{\oplus k_3}$ contribute
$2k_2+2k_3$ to $\dim W^H$ and $k_2+k_3$ to $\dim W^{H_1}$, while all other pure
$\TL(2,5)$-summands occurring in case $(1)$ contribute equally to $\dim W^H$ and
$\dim W^{H_1}$. The summands involving $\mathbb C_k$ contribute equally to
$\dim W^H$ and $\dim W^{H_1}$, because the only difference between $H$ and $H_1$
is in the $\TL(2,5)$-factor, where the corresponding fixed point dimensions are equal for
the summands under consideration. They contribute zero to $\dim W^{H_2}$, since $H_2$
contains the whole cyclic factor $C_k$. Finally, the summand
$W'\otimes \mathbb R_{\pm}$ contributes zero to all three fixed point dimensions, since
each of the subgroups $H,H_1,H_2$ contains the order two subgroup of the cyclic factor
on which $\mathbb R_{\pm}$ is non-trivial.

We split the proof into cases.

If
$$
        \dim W^H=0,
$$
then also
$$
        \dim W^{H_1}=\dim W^{H_2}=0.
$$
The subgroup $H={Q_8}_{(15)}\times C_{2^\ell}$ is a $2$-group, hence it satisfies the
group-theoretic hypothesis of \Cref{lemma:equal_dims}. Moreover, $H\leq H_1\cap H_2$, the subgroups
$H_1$ and $H_2$ generate $\TL(2,5)\times C_k$, and $\dim W^H=0$. Finally, both
$H_1=Q_{16}\times C_{2^\ell}$ and $H_2=(C_3\rtimes Q_8)\times C_k$ are non-Oliver groups, hence
they cannot admit one-fixed-point actions on homology spheres. Thus \Cref{lemma:equal_dims} gives a
contradiction.

Suppose next that
$$
        k_2+k_3=0
        \qquad\text{and}\qquad
        \dim W^H>0.
$$
Then
$$
        \dim W^H=\dim W^{H_1}>0
        \qquad\text{and}\qquad
        \dim W^{H_2}=0.
$$
Since $H$ is a $2$-group, Smith theory implies that $\Sigma^H$ is a positive-dimensional
$\mathbb Z_2$-homology sphere, in particular it is connected. Therefore Theorem 2.3
applies to
$$
        H={Q_8}_{(15)}\times C_{2^\ell},\qquad
        H_1=Q_{16}\times C_{2^\ell},\qquad
        H_2=(C_3\rtimes Q_8)\times C_k.
$$
Indeed, $H\leq H_1\cap H_2$, the subgroups $H_1$ and $H_2$ generate the whole group,
and $\dim W^H=\dim W^{H_1}$. Since $H_2$ is a non-Oliver group, it cannot act with exactly one
fixed point on a homology sphere. Thus \Cref{theorem:equaldims_Kamil} gives a contradiction.

Finally, suppose that
$$
        k_2+k_3>0.
$$
Then
$$
        \dim W^{H_2}=k_2+k_3>0.
$$
Also, the summands $V_{4,2}^{\oplus k_2}$ and $V_{5,1}^{\oplus k_3}$ contribute
$k_2+k_3$ to $\dim W^{H_1}$, hence
$$
        \dim W^{H_1}>0.
$$
Moreover,
$$
        \dim W^H=\dim W^{H_1}+\dim W^{H_2}.
$$
Since $H$ is a $2$-group, we may apply \Cref{theorem:pdis} with coefficients in $\mathbb Z_2$.
The subgroup inclusions and the generation condition were checked above. Hence
\Cref{theorem:pdis} gives a contradiction.

We now consider case $(2)$. We use the subgroups
$$
        H=C_8\times C_{2^\ell},\qquad
        H_1=Q_{16}\times C_{2^\ell},\qquad
        H_2=(C_5\rtimes C_8)\times C_k.
$$
Then
$$
        H\leq H_1\cap H_2
$$
and
$$
        \langle H_1,H_2\rangle=\TL(2,5)\times C_k.
$$
Using \Cref{table:dims_tl25} and \Cref{corollary:fixed_dims2}, we have
$$
        \dim W^{H_2}=k_3.
$$
Moreover,
$$
        \dim W^H=\dim W^{H_1}+k_3.
$$
Indeed, the summand $V_{5,2}^{\oplus k_3}$ contributes $2k_3$ to $\dim W^H$ and
$k_3$ to $\dim W^{H_1}$, while all remaining pure $\TL(2,5)$-summands occurring in
case $(2)$ contribute equally to $\dim W^H$ and $\dim W^{H_1}$. The summands involving
$\mathbb C_k$ contribute equally to $\dim W^H$ and $\dim W^{H_1}$ and contribute zero
to $\dim W^{H_2}$, because $H_2$ contains the whole cyclic factor $C_k$. The summand
$W'\otimes \mathbb R_{\pm}$ again contributes zero to all three fixed point dimensions,
since all three subgroups contain the order two subgroup of the cyclic factor on which
$\mathbb R_{\pm}$ is non-trivial.

Again we split the proof into cases.

If
$$
        \dim W^H=0,
$$
then also
$$
        \dim W^{H_1}=\dim W^{H_2}=0.
$$
The subgroup $H=C_8\times C_{2^\ell}$ is a $2$-group. Moreover,
$H\leq H_1\cap H_2$, the subgroups $H_1$ and $H_2$ generate $\TL(2,5)\times C_k$,
and $\dim W^H=0$. Since both $H_1=Q_{16}\times C_{2^\ell}$ and
$H_2=(C_5\rtimes C_8)\times C_k$ are non-Oliver groups, they cannot admit one-fixed-point actions
on homology spheres. Hence \Cref{lemma:equal_dims} gives a contradiction.

Suppose next that
$$
        k_3=0
        \qquad\text{and}\qquad
        \dim W^H>0.
$$
Then
$$
        \dim W^H=\dim W^{H_1}>0
        \qquad\text{and}\qquad
        \dim W^{H_2}=0.
$$
Since $H$ is a $2$-group, Smith theory implies that $\Sigma^H$ is a connected
positive-dimensional $\mathbb Z_2$-homology sphere. Therefore \Cref{theorem:equaldims_Kamil} applies to
$$
        H=C_8\times C_{2^\ell},\qquad
        H_1=Q_{16}\times C_{2^\ell},\qquad
        H_2=(C_5\rtimes C_8)\times C_k,
$$
yielding a contradiction.

Finally, suppose that
$$
        k_3>0.
$$
Then
$$
        \dim W^{H_2}=k_3>0.
$$
Also, the summand $V_{5,2}^{\oplus k_3}$ contributes $k_3$ to $\dim W^{H_1}$, hence
$$
        \dim W^{H_1}>0.
$$
Moreover,
$$
        \dim W^H=\dim W^{H_1}+\dim W^{H_2}.
$$
Since $H$ is a $2$-group, \Cref{theorem:pdis} applies with coefficients in $\mathbb Z_2$.
The subgroup inclusions and the generation condition were checked above, so \Cref{theorem:pdis}
gives a contradiction.

In all cases we obtain a contradiction. Therefore $W$ cannot occur as the tangential
representation at the fixed point of a one-fixed-point action on a homology sphere.
\end{proof}

\begin{theorem}\label{theorem:2main2_even}\emph{(cf. \Cref{theorem:2main})}
$\TL(2,5)\times C_k$ does not admit a one fixed point action on a homology sphere $\Sigma$ of dimension $n$ when $n\in\{6,\ldots,17\}\cup\{21\}$.
  
\end{theorem}
\begin{proof}
Analogously to the proof of \Cref{theorem:2main2}, combining both cases of \Cref{lemma:structural_tl25_even} and \Cref{lemma:one_dim}, we can conclude that $V$ must contain irreducible summands isomorphic to summands from each of the sets
$$
\mathcal{S}_1=\{V_{6},V_{4,2}\otimes\mathbb{C}_k,V_{5,1}\otimes\mathbb{C}_k,V_{6}\otimes \mathbb{C}_k\},
$$
and
$$
\mathcal{S}_2=\{V_{4,1},V_{6},V_{4,1}\otimes\mathbb{C}_k,V_{5,2}\otimes\mathbb{C}_k,V_{6}\otimes \mathbb{C}_k\}.
$$
Since the restriction of $V$ to $\TL(2,5)$ is effective as well, $V$ must also contain a summand from the set 
$$
\mathcal{F}=\{V_{8,i},V_{8,j}\otimes_{\mathbb{C}}\mathbb{C}_k,V_{8,l}\otimes\mathbb{R}_{\pm},V_{12},V_{12}\otimes_{\mathbb{C}}\mathbb{C}_k,V_{12}\otimes\mathbb{R}_{\pm}\},\quad i,j,l\in\{1,2\}.
$$
The above restrictions already exclude the cases $n\in\{6,\ldots,13\}\cup\{15,17\}$. We are therefore left with the cases $n=14$, $n=16$, and $n=21$. 

Let $n=14$. We note that any $W\in\mathcal{F}$ contains the diagonally embedded $C_2$ in its kernel. The only possibility to get dimension $14$ and satisfy the constraints mentioned above is to add $V_6$ to $W$. The representation $V_6$ contains the same diagonally embedded $C_2$ in its kernel, and thus we cannot obtain a faithful representation $V$ in this case.

For $n\in\{16,21\}$ we group all the possible representations similarly to the proof of \Cref{theorem:2main2}. The additional representations we have to consider are obtained from those from the proof of \Cref{theorem:2main2} by substituting the $5$ and $8$-dimensional direct summands with $V_{5,i}\otimes\mathbb{R}_{\pm}$ and $V_{8,i}\otimes\mathbb{R}_{\pm}$ respectively. Arguing similarly to the proof of \Cref{theorem:2main2}, we consider $W$, the tangent space to $\Sigma'=\Sigma^{C_2}$ at the fixed point, where $C_2\trianglelefteq\TL(2,5)$ (note that $W$ is an $S_5\times C_k$-representation). This leaves us with the cases
$$
W=V_6\oplus(V_{5,i}\otimes\mathbb{R}_{\pm})\oplus\mathbb{C}_k\quad\text{and}\quad W=V_6\oplus(V_{5,i}\otimes\mathbb{R}_{\pm})\oplus(V_1\otimes\mathbb{C}_k).
$$
Let $C_{k/2}$ be the cyclic subgroup of $C_k$ of index $2$. By considering the intersection number $\mathrm{mod}$ $2$ of $(\Sigma')^{A_5}$ and $(\Sigma')^{C_{k/2}}$ in $\Sigma'$, we conclude that $(\Sigma')^{A_5\times C_{k/2}}$ consists of an even number of points. The assertion follows then from \Cref{lemma:index_two_excl2} applied with parameters $G=S_5\times C_k$, $H=A_5\times C_{k/2}$, $H_1=A_5$, and $H_2=C_{k/2}$.

\end{proof}

\noindent Similarly to before (see \Cref{lemma:general_form_sl25_even} and \Cref{lemma:structural_tl25_even}), we get the following.
\begin{lemma}\label{lemma:structural_s5_even}
Let $W$ be one of the following representations of $S_{5}\times C_k$:
\begin{enumerate}
    \item $   V_{4,2}^{\oplus k_{1}}   \oplus V_{5,1}^{\oplus k_{2}}   \oplus V_{5,2}^{\oplus k_{3}} \oplus\mathbb{C}_k^{\oplus k_{4}} \oplus (V_{1}\otimes \mathbb{C}_k)^{\oplus k_{5}} \oplus (V_{4,2}\otimes \mathbb{C}_k)^{\oplus k_{6}}\oplus (V_{5,1}\otimes \mathbb{C}_k)^{\oplus k_{7}}\oplus (W'\otimes\mathbb{R}_{\pm})  $,
 
    \item  $V_{4,1}^{\oplus k_{1}}   \oplus V_{4,2}^{\oplus k_{2}} \oplus V_{5,1}^{\oplus k_{3}}\oplus V_{5,2}^{\oplus k_{4}}  \oplus\mathbb{C}_k^{\oplus k_{5}} \oplus (V_{1}\otimes \mathbb{C}_k)^{\oplus k_{6}}\oplus (V_{4,1}\otimes \mathbb{C}_k)^{\oplus k_{7}}\oplus (V_{5,2}\otimes \mathbb{C}_k)^{\oplus k_{8}}\oplus (W'\otimes\mathbb{R}_{\pm})     $,

\end{enumerate}
where $W'$ is a real $S_5$-representation. Then $W$ cannot occur as the tangential
representation at the fixed point of a one fixed point action
on a homology sphere where $k_{i}$
are arbitrary non-negative integers.

\end{lemma}
\begin{proof}
Assume, for a contradiction, that there exists a homology sphere $\Sigma$ with a
one-fixed-point action of $S_5\times C_k$ and with tangential representation at the fixed
point isomorphic to $W$. Write
$$
        C_k\cong C_{2^\ell}\times C_m,
$$
where $m$ is odd. Since $k$ is even, we have $\ell\geq 1$.

We first consider case $(1)$. We use the subgroups
$$
        H=C_4\times C_{2^\ell},\qquad
        H_1=D_8\times C_{2^\ell},\qquad
        H_2=(C_5\rtimes C_4)\times C_k.
$$
Then
$$
        H\leq H_1\cap H_2
$$
and
$$
        \langle H_1,H_2\rangle=S_5\times C_k.
$$
Using \Cref{table:dims_s5} and \Cref{corollary:fixed_dims2}, we have
$$
        \dim W^{H_2}=k_3.
$$
Moreover,
$$
        \dim W^H=\dim W^{H_1}+k_3.
$$
Indeed, the summand $V_{5,2}^{\oplus k_3}$ contributes $2k_3$ to $\dim W^H$ and
$k_3$ to $\dim W^{H_1}$, while all other pure $S_5$-summands occurring in case $(1)$
contribute equally to $\dim W^H$ and $\dim W^{H_1}$. The summands involving
$\mathbb C_k$ also contribute equally to $\dim W^H$ and $\dim W^{H_1}$, because the
only difference between $H$ and $H_1$ is in the $S_5$-factor, where the corresponding
fixed point dimensions are equal for the summands under consideration. They contribute
zero to $\dim W^{H_2}$, since $H_2$ contains the whole factor $C_k$. Finally, the summand
$W'\otimes \mathbb R_{\pm}$ contributes zero to all three fixed point dimensions, since
each of the subgroups $H,H_1,H_2$ contains the order two subgroup of the cyclic factor
on which $\mathbb R_{\pm}$ is non-trivial.

We now split the proof into cases.

If
$$
        \dim W^H=0,
$$
then also
$$
        \dim W^{H_1}=\dim W^{H_2}=0.
$$
The subgroup $H=C_4\times C_{2^\ell}$ is a $2$-group, hence it satisfies the
group-theoretic hypothesis of \Cref{lemma:equal_dims}. Moreover, $H\leq H_1\cap H_2$, the subgroups
$H_1$ and $H_2$ generate $S_5\times C_k$, and $\dim W^H=0$. Finally, both
$H_1=D_8\times C_{2^\ell}$ and $H_2=(C_5\rtimes C_4)\times C_k$ are non-Oliver groups, hence
they cannot admit one-fixed-point actions on homology spheres. Thus \Cref{lemma:equal_dims} gives a
contradiction.

Suppose next that
$$
        k_3=0
        \qquad\text{and}\qquad
        \dim W^H>0.
$$
Then
$$
        \dim W^H=\dim W^{H_1}>0
        \qquad\text{and}\qquad
        \dim W^{H_2}=0.
$$
Since $H$ is a $2$-group, Smith theory implies that $\Sigma^H$ is a positive-dimensional
$\mathbb Z_2$-homology sphere, in particular it is connected. Therefore \Cref{theorem:equaldims_Kamil}
applies to
$$
        H=C_4\times C_{2^\ell},\qquad
        H_1=D_8\times C_{2^\ell},\qquad
        H_2=(C_5\rtimes C_4)\times C_k.
$$
Indeed, $H\leq H_1\cap H_2$, the subgroups $H_1$ and $H_2$ generate the whole group,
and $\dim W^H=\dim W^{H_1}$. Since $H_2$ is non-Oliver, it cannot act with exactly one
fixed point on a homology sphere. Thus \Cref{theorem:equaldims_Kamil} gives a contradiction.

Finally, suppose that
$$
        k_3>0.
$$
Then
$$
        \dim W^{H_2}=k_3>0
$$
and, since $V_{5,2}^{\oplus k_3}$ contributes $k_3$ to $\dim W^{H_1}$, we also have
$$
        \dim W^{H_1}>0.
$$
Moreover,
$$
        \dim W^H=\dim W^{H_1}+\dim W^{H_2}.
$$
Since $H$ is a $2$-group, we may apply \Cref{theorem:pdis} with coefficients in $\mathbb Z_2$.
The subgroup inclusions and the generation condition were checked above. Hence
\Cref{theorem:pdis} gives a contradiction.

We now consider case $(2)$. We use the subgroups
$$
        H={C_2^2}_{(15)}\times C_{2^\ell},\qquad
        H_1=D_8\times C_{2^\ell},\qquad
        H_2=D_{12}\times C_k.
$$
Then
$$
        H\leq H_1\cap H_2
$$
and
$$
        \langle H_1,H_2\rangle=S_5\times C_k.
$$
Using \Cref{table:dims_s5} and \Cref{corollary:fixed_dims2}, we have
$$
        \dim W^{H_2}=k_2+k_3.
$$
Moreover,
$$
        \dim W^H=\dim W^{H_1}+k_2+k_3.
$$
Indeed, the summands $V_{4,2}^{\oplus k_2}$ and $V_{5,1}^{\oplus k_3}$ contribute
$2k_2+2k_3$ to $\dim W^H$ and $k_2+k_3$ to $\dim W^{H_1}$, while all remaining
pure $S_5$-summands occurring in case $(2)$ contribute equally to $\dim W^H$ and
$\dim W^{H_1}$. The summands involving $\mathbb C_k$ contribute equally to
$\dim W^H$ and $\dim W^{H_1}$ and contribute zero to $\dim W^{H_2}$, because
$H_2$ contains the whole cyclic factor $C_k$. The summand $W'\otimes \mathbb R_{\pm}$
again contributes zero to all three fixed point dimensions, since all three subgroups contain
the order two subgroup of the cyclic factor on which $\mathbb R_{\pm}$ is non-trivial.

Again we split the proof into cases.

If
$$
        \dim W^H=0,
$$
then also
$$
        \dim W^{H_1}=\dim W^{H_2}=0.
$$
The subgroup $H={C_2^2}_{(15)}\times C_{2^\ell}$ is a $2$-group. Moreover,
$H\leq H_1\cap H_2$, the subgroups $H_1$ and $H_2$ generate $S_5\times C_k$, and
$\dim W^H=0$. Since both $H_1=D_8\times C_{2^\ell}$ and $H_2=D_{12}\times C_k$
are non-Oliver, they cannot admit one-fixed-point actions on homology spheres. Hence
\Cref{lemma:equal_dims} gives a contradiction.

Suppose next that
$$
        k_2+k_3=0
        \qquad\text{and}\qquad
        \dim W^H>0.
$$
Then
$$
        \dim W^H=\dim W^{H_1}>0
        \qquad\text{and}\qquad
        \dim W^{H_2}=0.
$$
Since $H$ is a $2$-group, Smith theory implies that $\Sigma^H$ is a connected
positive-dimensional $\mathbb Z_2$-homology sphere. Therefore \Cref{theorem:equaldims_Kamil} applies to
$$
        H={C_2^2}_{(15)}\times C_{2^\ell},\qquad
        H_1=D_8\times C_{2^\ell},\qquad
        H_2=D_{12}\times C_k,
$$
yielding a contradiction.

Finally, suppose that
$$
        k_2+k_3>0.
$$
Then
$$
        \dim W^{H_2}=k_2+k_3>0.
$$
Also, the summands $V_{4,2}^{\oplus k_2}$ and $V_{5,1}^{\oplus k_3}$ contribute
$k_2+k_3$ to $\dim W^{H_1}$, hence
$$
        \dim W^{H_1}>0.
$$
Moreover,
$$
        \dim W^H=\dim W^{H_1}+\dim W^{H_2}.
$$
Since $H$ is a $2$-group, \Cref{theorem:pdis} applies with coefficients in $\mathbb Z_2$.
The subgroup inclusions and the generation condition were checked above, so \Cref{theorem:pdis}
gives a contradiction.

In all cases we obtain a contradiction. Therefore $W$ cannot occur as the tangential
representation at the fixed point of a one-fixed-point action on a homology sphere.
\end{proof}

\begin{theorem}\label{theorem:3main22}\emph{(cf. \Cref{theorem:3main})}
$S_{5}\times C_k$ does not admit a one fixed point action on a homology sphere $\Sigma$ of dimension $n$ where $n\in\{6,\ldots,13\}\cup\{15,17\}$.  
\end{theorem}

\begin{proof}
Let $V=T_x\Sigma$, where $x$ is the only fixed point of the action of $S_5\times C_k$ on $\Sigma$ (we reason, as usual, by contradiction). Similarly as in the proofs of \Cref{theorem:3main2}, combining both cases of \Cref{lemma:structural_s5_even} and \Cref{lemma:one_dim}, we can conclude that $V$ must contain irreducible summands isomorphic to summands from each of the sets
$$
\mathcal{S}_1=\{V_{4,1},V_6,V_{4,1}\otimes\mathbb{C}_k,V_{5,2}\otimes\mathbb{C}_k,V_{6}\otimes \mathbb{C}_k\},
$$
and
$$
\mathcal{S}_2=\{V_6,V_{4,2}\otimes\mathbb{C}_k,V_{5,1}\otimes\mathbb{C}_k,V_{6}\otimes \mathbb{C}_k\}.
$$

\noindent We explain the reduction to the dimensions $12,13,15$, and $17$. The restriction that $V$ contains a summand from both $\mathcal{S}_1$ and $\mathcal{S}_2$ enforces $V$ to be of dimension at least $6$ (we can fulfill this requirement by $V$ having $V_6$ as a direct summand). Moreover, by \Cref{table:dims_s5} and the fact that $V$ must not contain an irreducible summand of dimension $1$, the smallest possible dimension of another summand of $V$ is $4$. This reasoning rules out all dimensions up to $9$. It is also not hard to see that dimension $10$ is not possible: the only configuration potentially possible would be $V=V_6\oplus W$, where $W$ is some irreducible $4$-dimensional $S_5\times C_k$-representation. By \Cref{table:dims_s5} this leaves, however, the non-effective cases only. Analogously, we can rule out dimension $11$ as in that case $V$ must have contained an irreducible summand of dimension $5$. 



 In dimension $12$ we are only left with $V\in\{V_{4,1}\oplus(V_{4,2}\otimes\mathbb{C}_k),V_6\oplus V_{4,i}\oplus\mathbb{C}_k,V_6\oplus(V_{4,i}\otimes\mathbb{R}_{\pm})\oplus\mathbb{C}_k\}$. The cases with $\mathbb{C}_k$ as a direct summand can be excluded by \Cref{corolarry:one_dims_struct}. In the case $V=V_{4,1}\oplus(V_{4,2}\otimes\mathbb{C}_k)$ the argument from the proof of \Cref{theorem:3main2} for $V_{4,1}\oplus(V_{4,2}\otimes\mathbb{C}_k)$ can be applied verbatim.

 In dimension $13$ the only possibilities are $V=V_6\oplus (V_{5,i}\otimes\mathbb{R}_{\pm})\oplus\mathbb{C}_k$ and $V=V_6\oplus (V_{5,i}\otimes\mathbb{R}_{\pm})\oplus(V_1\otimes\mathbb{C}_k)$. These cases follow by \Cref{lemma:index_two_excl2} applied with parameters $G=S_5\times C_k$, $H=A_5\times C_{k/2}$, $H_1=A_5$, and $H_2=C_{k/2}$. 

In dimension $15$, the remaining candidates are obtained from the candidates in
dimension $13$ by adding one more two-dimensional summand, isomorphic either to
$\mathbb C_k$ or to $V_1\otimes \mathbb C_k$. Therefore the possible representations are
$$
        V\cong V_6\oplus (V_{5,i}\otimes \mathbb R_{\pm})
        \oplus \mathbb C_k^{\oplus 2},
$$
$$
        V\cong V_6\oplus (V_{5,i}\otimes \mathbb R_{\pm})
        \oplus \mathbb C_k\oplus (V_1\otimes \mathbb C_k),
$$
and
$$
        V\cong V_6\oplus (V_{5,i}\otimes \mathbb R_{\pm})
        \oplus (V_1\otimes \mathbb C_k)^{\oplus 2},
$$
where $i=1,2$. These cases are excluded by the same argument as in dimension $13$. Indeed, for each of
the three representations above we apply \Cref{lemma:index_two_excl2} to the subgroups
$$
        H=A_5\times C_{k/2},\qquad
        H_1=A_5,\qquad
        H_2=C_{k/2}.
$$
Then
$$
        \langle H_1,H_2\rangle=H,\qquad
        H\triangleleft S_5\times C_k,
$$
and
$$
        (S_5\times C_k)/H\cong C_2\times C_2.
$$
Moreover, in each of the three cases we have
$$
        \dim V^{A_5}=4
        \qquad\text{and}\qquad
        \dim V^{C_{k/2}}=11.
$$
Thus
$$
        \dim V^{A_5}+\dim V^{C_{k/2}}=15=\dim V.
$$
Therefore \Cref{lemma:index_two_excl2} applies and excludes all the dimension $15$ candidates listed above.

 For $17$-dimensional representation $V$ such that 
 $\mathrm{Res}_{S_5 \times C_{k/{2^l}}}^{S_5 \times C_{k}}(V)\cong  V_{5,i}\oplus V_{4,1}\oplus(V_{4,2}\otimes\mathbb{C}_k)$ by \Cref{lemma:37} we know that $\Sigma^{S_5 \times C_{k/{2^l}}}$, where $k=2^lm$ for odd $m$, contains an even number of points, so the same remains true for $\Sigma^{S_5 \times C_{k}}= (\Sigma^{S_5 \times C_{k/{2^l}}})^{C_{2^l}}$. Applying \Cref{corolarry:one_dims_struct} we can restrict ourselves to the cases $V=(V_{6}\otimes\mathbb{C}_{k})\oplus V_{5,i}$ and $V=(V_{6}\otimes\mathbb{C}_{k})\oplus (V_{5,i}\otimes \mathbb{R}_{\pm})$.

 Let us now assume that $V=(V_{6}\otimes\mathbb{C}_{k})\oplus V_{5,i}$ for $i=1,2$. To avoid confusion we denote by $\mathcal{C}_{2}$ the subgroup of $C_{k}$ of order two. Consider the restriction of $V$ to $A_5\times \mathcal{C}_{2}$-representation $W=\mathrm{Res}^{S_5 \times C_{k}}_{A_5 \times \mathcal{C}_{2}}(V) \cong V_5 \oplus (V_{3}\otimes \mathbb{R}_{\pm})^{\oplus 4}$ (see \Cref{table:fpsdim_A5}). By inspecting the table of fixed points we have $\mathrm{dim}W^{C_{2} \times \mathcal{C}_{2}}=3, \mathrm{dim}W^{C_{2}^{2} \times \mathcal{C}_{2}}=2$, and $\mathrm{dim}W^{S_3 \times \mathcal{C}_{2}}=1$ which by the $\mathrm{mod}$ $2$ intersection number argument shows that $\Sigma^{A_5 \times \mathcal{C}_{2}}$ is not a single point, which together with $\Sigma^{A_4}=S^{0}$ (see the proof of \Cref{lemma:3_8}) shows that $\Sigma^{A_5 \times \mathcal{C}_{2}} = S^{0}$. Finally $\Sigma^{S_5 \times C_k}=(\Sigma^{A_5 \times \mathcal{C}_{2}})^{C_{2} \times C_{k/2}} = S^{0}$.


 Now consider
$
        V=(V_6\otimes \mathbb C_k)\oplus (V_{5,i}\otimes \mathbb R_{\pm})$, $i=1,2$.
We may assume that the kernel of the $C_k$-representation occurring in the summand
$V_6\otimes \mathbb C_k$ is the subgroup $\mathcal{C}_2$. Otherwise, a similar
argument as in the previous case applies ($\mathrm{dim}V^{C_{2} \times \mathcal{C}_{2}}=\mathrm{dim}V^{C_{2}^{2} \times \mathcal{C}_{2}}=\mathrm{dim}V^{S_3 \times \mathcal{C}_{2}}=0$).

Put
$$
        L={C_2^2}_{(15)}\times C_{2^l},\qquad
        K_1=D_8\times C_k,\qquad
        K_2=D_{12}\times C_k.
$$
Then
$$
        L\leq K_1\cap K_2
$$
and
$$
        \langle K_1,K_2\rangle=S_5\times C_k.
$$

\noindent Hence,
$$
        \Sigma^L
        =
        \Sigma^{{C_2^2}_{(15)}\times C_k}
        \cong S^0.
$$

\noindent Now recall that neither $D_8\times C_k$ nor $D_{12}\times C_k$ admits a one-fixed-point
action on a homology sphere. Since
$$
        \Sigma^{K_1}\subseteq \Sigma^L\cong S^0
        \qquad\text{and}\qquad
        \Sigma^{K_2}\subseteq \Sigma^L\cong S^0,
$$
and both fixed point sets contain the point $x$, each of them is either $\{x\}$ or the whole
of $S^0$. The first possibility is excluded by the non-existence of one-fixed-point actions for
$K_1$ and $K_2$. Hence,
$$
        \Sigma^{K_1}=\Sigma^L\cong S^0
        \qquad\text{and}\qquad
        \Sigma^{K_2}=\Sigma^L\cong S^0.
$$
Therefore
$$
        \Sigma^{S_5\times C_k}
        =
        \Sigma^{\langle K_1,K_2\rangle}
        =
        \Sigma^{K_1}\cap \Sigma^{K_2}
        =
        \Sigma^L
        \cong S^0.
$$
In particular, $\Sigma^{S_5\times C_k}$ is not a singleton, which contradicts the assumption
that the action has exactly one fixed point.

    \end{proof}

\section{Appendix -- the case of $A_5\times C_k$}
One of the main goals of this section is to give a correct proof of \cite[Proposition 7.1]{Morimoto2025_OneFixedPointIII}, case $n=10$ for the partition $p(V)=(4,6)$. The proof given in \cite[Proposition 7.1]{Morimoto2025_OneFixedPointIII} uses the equality $\Sigma^{A_5}=\Sigma^{C_2^2}$ which does not hold in the considered case.
\begin{lemma}\label{theorem:a5_thm}
    Let $A_{5}$ act on a $4$-dimensional $\mathbb{Z}_{p}$-homology sphere $\Sigma$, where $p$ is an arbitrary prime, with a nonempty, discrete set of fixed points. Then $\Sigma^{A_{5}}$ contains $kp + 2m$ points for some nonnegative integers $k$ and $m$.

\end{lemma}

\begin{proof}
First assume that $p$ is odd. Recall that $A_5$ is generated by two subgroups of order $3$, for instance by the subgroup $C_3$ generated by $(1,2,3)$ and the subgroup $\mathcal{C}_3$ generated by $(3,4,5)$. Both $\Sigma^{C_3}$ and $\Sigma^{\mathcal{C}_3}$ are disjoint unions of orientable surfaces and, possibly, isolated points, since they are fixed point sets of actions of odd-order groups on an orientable manifold. In our argument, these isolated points may be ignored, because they cannot be fixed by the whole group. Indeed, if an isolated point were $A_5$-fixed, then the tangent representation $V$ of $A_5$ at this point would satisfy $\dim V^{C_3}=0$; however, no four-dimensional real representation of $A_5$ with this property exists.

The cohomological intersection number modulo $p$ of $\Sigma^{C_3}$ and $\Sigma^{\mathcal{C}_3}$ in $\Sigma$ is equal to $0$. Moreover, we restrict our attention to the case where the relevant intersections are discrete. By the interpretation of the cohomological intersection number, each path component of $\Sigma^{C_3}$ intersects each component of $\Sigma^{\mathcal{C}_3}$ in $pk$ points, for some integer $k$, together with possibly an even number of additional points, the latter arising from intersections which occur once with positive sign and once with negative sign. Counting over all path components of $\Sigma^{C_3}$, the result follows.

For $p=2$, the same argument applies, but without the need to discuss orientations of submanifolds or signs of intersections.
\end{proof}

\begin{table}[H]
\begin{center}
{\fontsize{10pt}{12pt}\selectfont 
\begin{tabular}{|c||c|c|c|c|c|c|c|c|c|c|} \hline
$H$ & 
$E$ & $C_{2}$ & $C_{3}$ &$C_{5}$ & 
$(C_{2})^2$ & $S_{3}$ & $D_{10}$ &$A_{4}$ & 
$A_{5}$    \\ 
\hline
\hline

$V_{3}$ & 
$3$ & $1$ & $1$ & $1$ & $0$ & 
$0$ & $0$ & $0$ & $0$ \\
\hline 
$V_{4}$ & 
$4$ & $2$ & $2$ & $0$ & $1$ & 
$1$ & $0$ & $1$ & $0$  \\
\hline

$V_{5}$ & 
$5$ & $3$ & $1$ & $1$ & $2$ & 
$1$ & $1$ & $0$ & $0$   \\
\hline

\end{tabular}
}
\vskip 0.2cm
\caption{Fixed point dimension table for the non-trivial real irreducible representations of $A_5$.}\label{table:fpsdim_A5}

\end{center}
\end{table}

 \begin{corollary}\label{corolarry:a5even}
Let $\Sigma$ be a $4$-dimensional $\mathbb{Z}_{p}$-homology sphere, where $p$ is an arbitrary prime. If $A_{5}$ acts on $\Sigma$ as in \Cref{theorem:a5_thm}, then $\Sigma^{A_{5}}$ contains an even number of points.
\end{corollary}

\begin{proof}
Since $\Sigma$ is a compact manifold, its integral homology groups are finitely generated. Therefore only finitely many primes can occur as torsion primes in the homology of $\Sigma$. It follows that there exists an integer $q$ such that $\Sigma$ is a $\mathbb{Z}_{r}$-homology sphere for every prime $r>q$. The assertion now follows from \Cref{theorem:a5_thm}.
\end{proof}

\begin{theorem}\label{theorem:a5ck_exclusion}
Let $k>2$. Then the group $A_{5}\times C_{k}$ does not admit a one-fixed-point action on an
$n$-dimensional homology sphere for $n\in\{7,8,10\}$.
\end{theorem}

\begin{proof}
Assume, for a contradiction, that $A_{5}\times C_{k}$ acts on a homology sphere $\Sigma$
with exactly one fixed point $x$, and let
$$
        V = T_x\Sigma
$$
be the tangential $A_{5}\times C_{k}$-representation at $x$.

First suppose that $k$ is even. By \Cref{lemma:index_two_excl2}, applied to the index-two
subgroup
$$
        A_5\times C_{k/2}\triangleleft A_5\times C_k
$$
and to the subgroups $H_1=C_{k/2}$ and $H_2=A_5$, we can exclude every representation
of the form
$$
        W_1\oplus W_2,
$$
where $W_1\cong \mathbb{C}_k^{m}$ for some positive integer $m$, and $W_2$ is a non-trivial direct
sum of irreducible $A_{5}\times C_{2}$-representations. If $k$ is odd, the analogous reduction
is obtained by omitting the $C_2$-factor. The case when either of $W_i$ equals $\{0\}$ can be excluded by the effectiveness assumption.

Moreover, since the action has exactly one fixed point, the tangential representation $V$
contains no trivial summand. By \Cref{lemma:one_dim}, $V$ contains no non-trivial one-dimensional
irreducible summand either. Since $V$ must be faithful, we may also exclude all representations whose irreducible summands are lifted from the
projection
$
        A_5\times C_k \twoheadrightarrow A_5.
$
After these reductions, the only remaining possibilities are
$$
        V_4\otimes \mathbb{C}_k,\qquad
        V_5\otimes \mathbb{C}_k,\qquad
        (V_3\otimes \mathbb{C}_k)\oplus V_4,
$$
and, when $k$ is even,
$$
        (V_3\otimes \mathbb{C}_k)\oplus (V_4\otimes \mathbb{R}_{\pm}).
$$

Consider first the cases
$$
        V\cong V_4\otimes \mathbb{C}_k
        \quad\text{or}\quad
        V\cong V_5\otimes \mathbb{C}_k .
$$
Choose a subgroup $C_{p^l}\leq C_k$ of prime power order such that
$$
        \dim V^{C_{p^l}}=\dim V^{C_k}=0.
$$
By Smith theory, $\Sigma^{C_{p^l}}$ is a $0$-dimensional $\mathbb{Z}_p$-homology sphere, hence
$$
        \Sigma^{C_{p^l}}\cong S^0.
$$
Since $C_k/C_{p^l}$ acts on this two-point set and fixes $x$, it fixes both points. Therefore
$$
        \Sigma^{C_k}=\Sigma^{C_{p^l}}\cong S^0.
$$
It follows that
$$
        \Sigma^{A_5\times C_k}
        =
        \bigl(\Sigma^{C_k}\bigr)^{A_5}
        =
        (S^0)^{A_5}
        \neq \{x\},
$$
which contradicts the assumption that the action has exactly one fixed point.

Now consider
$$
        V\cong (V_3\otimes \mathbb{C}_k)\oplus V_4 .
$$
Again choose a subgroup $C_{p^l}\leq C_k$ of prime power order such that
$$
        \dim V^{C_{p^l}}=\dim V^{C_k}=4.
$$
Then $\Sigma^{C_{p^l}}$ is a $4$-dimensional $\mathbb{Z}_p$-homology sphere. Since
$C_k/C_{p^l}$ acts on $\Sigma^{C_{p^l}}$ and the fixed point set of this action has the same tangent
dimension at $x$ as $\Sigma^{C_{p^l}}$ itself, we get
$$
        \Sigma^{C_k}=\Sigma^{C_{p^l}}.
$$
Hence $\Sigma^{C_k}$ is a $4$-dimensional $\mathbb{Z}_p$-homology sphere with an induced
$A_5$-action. By \Cref{corolarry:a5even},
$$
        \bigl(\Sigma^{C_k}\bigr)^{A_5}
        =
        \Sigma^{A_5\times C_k}
$$
contains an even number of points. In particular, it cannot be a singleton.

It remains to treat the case
$$
        V\cong (V_3\otimes \mathbb{C}_k)\oplus (V_4\otimes \mathbb{R}_{\pm}),
$$
which can occur only when $k$ is even. Choose a subgroup $C_{p^l}\leq C_{k/2}$ of prime power order.
Then
$$
        \dim V^{C_{p^l}}=\dim V^{C_{k/2}}=4.
$$
As above, Smith theory and the dimension equality imply that
$$
        \Sigma^{C_{k/2}}=\Sigma^{C_{p^l}},
$$
so $\Sigma^{C_{k/2}}$ is a $4$-dimensional $\mathbb{Z}_p$-homology sphere. Applying
\Cref{corolarry:a5even} to the induced $A_5$-action on $\Sigma^{C_{k/2}}$, we conclude that
$$
        \Sigma^{A_5\times C_{k/2}}
$$
contains an even number of points.

Now set
$$
        H=A_5\times C_{k/2}.
$$
Then $H\triangleleft A_5\times C_k$ and
$$
        (A_5\times C_k)/H\cong C_2.
$$
The quotient $C_2$ acts on $\Sigma^{H}$, which contains an even number of points, and its fixed point set is precisely
$$
        \Sigma^{A_5\times C_k}.
$$
Since a $C_2$-action on a finite set has the number of fixed points congruent modulo $2$ to
the cardinality of the whole set, and since $x$ is fixed, the set
$\Sigma^{A_5\times C_k}$ cannot consist of exactly one point.

In every possible case we obtained a contradiction. Therefore $A_5\times C_k$ does not admit
a one-fixed-point action on an $n$-dimensional homology sphere for
$n\in\{7,8,10\}$.
\end{proof}

\section{Appendix -- actions of extensions of finite groups by a torus}
 In this appendix we give a restriction on the dimension of homology spheres admitting one fixed point actions of $G$, where $G$ is an extension of a finite Oliver group by the $m$-dimensional torus $T^m \cong (S^1)^m$.

The additional information which will be used is Smith theory for torus actions, namely: if $T^n$ acts on a $\mathbb{Z}$-homology sphere $\Sigma$, then $\Sigma^{T^n}$ is also a $\mathbb{Z}$-homology sphere.

For induction steps in proofs we will use a sequence of normal subgroups of $T^m$
\[
S^1 = T^1 \triangleleft T^2 \triangleleft \dots \triangleleft T^m.
\]






\begin{theorem}\label{theorem:51}
Let $H$ be a finite Oliver group which does not admit a one-fixed-point action (not necessarily effective) on a
$k$-dimensional homology sphere for $k<r$. If $G$ fits into an exact sequence
$$
1\rightarrow T^m \rightarrow G \rightarrow H\rightarrow 1,
$$
then $G$ does not admit a one-fixed-point action on an $n$-dimensional homology sphere
if $n<r+2m$.
\end{theorem}

\begin{proof}
Suppose, for a contradiction, that such an action exists. Let $\Sigma$ be an
$n$-dimensional homology sphere with a one-fixed-point $G$-action, and let $x$ be its
unique $G$-fixed point. Put
$$
        V=T_x\Sigma .
$$
Thus $V$ is the tangential real $G$-representation at $x$.

Denote by $T^m\triangleleft G$ the torus appearing in the exact sequence and $F=\Sigma^{T^m}$. 
We show first that
$$
    \dim V-\dim V^{T^m}\geq 2m.
$$
Indeed, consider the restriction of the tangential representation to the torus,
$$
        \operatorname{Res}^{G}_{T^m}V.
$$
If some non-trivial element of $T^m$ acted trivially on $V$, then, by the Slice Theorem, it
would act trivially on a neighbourhood of $x$. Since the action may be taken to be by
isometries with respect to a $G$-invariant Riemannian metric, this element would act
trivially on all of $\Sigma$, contradicting the effectiveness of the $G$-action. Hence
$\operatorname{Res}^{G}_{T^m}V$ is a faithful real representation of $T^m$.

Every non-trivial irreducible real representation of a torus has real dimension $2$.
Moreover, a faithful representation of $T^m$ must contain at least $m$ independent
non-trivial torus characters. Therefore the non-trivial part of
$\operatorname{Res}^{G}_{T^m}V$ has dimension at least $2m$, and so
$$
        \dim V^{T^m}\leq \dim V-2m=n-2m.
$$
Since $n<r+2m$, we obtain
$$
        \dim F=\dim V^{T^m}\leq n-2m<r.
$$
\noindent Now, because $T^m$ is normal in $G$, the action of $G$ on $F=\Sigma^{T^m}$ factors through
the quotient group
$$
        G/T^m\cong H.
$$
Furthermore,
$$
        F^H
        =
        \left(\Sigma^{T^m}\right)^H
        =
        \Sigma^G
        =
        \{x\}.
$$
Thus $H$ admits a one-fixed-point action on the homology sphere $F$, whose dimension is
strictly smaller than $r$. This contradicts the assumption on $H$.

Therefore $G$ does not admit a one-fixed-point action on an $n$-dimensional homology
sphere for $n<r+2m$.
\end{proof}

More can be said in the case of groups $T^m \times H$, where $H$ is fixed, for example $A_5$ or $S_5$. The next lemma is an extension of the previous results, where $T^m$ is replaced by a finite cyclic group $C_k$.

\begin{lemma}\label{lemma:52}
Let $V$ be a $(T^m \times H)$-representation isomorphic to $V \cong W_1 \oplus W_2$, where $W_1$ (respectively $W_2$) is a pure $T^m$- (respectively $H$)-representation. Then $V$ cannot occur as the tangential representation of a one fixed point action on a homology sphere.

\end{lemma}

\begin{proof}
Let $\dim V = n$. Then
\[
\dim V^{T^m} + \dim V^H = n,
\]
and both summands have positive dimension by the effectiveness assumption. Therefore, the result follows by standard intersection number arguments.

\end{proof}

\begin{theorem}\label{theorem:53}
Let $H$ be a finite Oliver group such that:
\begin{enumerate}
    \item no non-trivial irreducible real $H$-representation admits an invariant complex
    structure; equivalently, if $\rho\colon H\rightarrow O(\dim W)$ is such a representation and
    $\dim W$ is even, then the image of $\rho$ is not conjugate into $U(\dim W/2)$;

    \item $H$ has no non-trivial irreducible real representations of dimension $d$ with
    $1<d<k$;

    \item $H$ does not admit a one-fixed-point action (not necessarily effective) on any homology sphere of dimension
    less than $m$.
\end{enumerate}
Then $T^n\times H$ does not admit a one-fixed-point action on a homology sphere of
dimension less than $m+2nk$.
\end{theorem}

\begin{proof}
We argue by induction on $n$. First consider the case $n=1$, so that the group is
$S^1\times H$. Suppose that $S^1\times H$ acts with exactly one fixed
point $x$ on a homology sphere $\Sigma$, and put
$$
        V=T_x\Sigma .
$$
Thus $V$ is the tangential real $S^1\times H$-representation at $x$.

Let
$
        F=\Sigma^{S^1}
$. Then $F$ is a
homology sphere and
$$
        T_xF=V^{S^1}.
$$
Moreover, the induced $H$-action on $F$ has
$$
        F^H=\Sigma^{S^1\times H}=\{x\}.
$$
By the assumption on $H$, we must therefore have
$$
        \dim V^{S^1}=\dim F\geq m.
$$

\noindent Now decompose the $S^1\times H$-representation $V$ as
$$
        V=V^{S^1}\oplus N,
$$
where $N$ is the sum of the irreducible summands on which $S^1$ acts non-trivially.
Since the original $S^1\times H$-action is effective, we have $N\neq 0$.

By the representation-theoretic assumptions on $H$ and \Cref{lemma:52} every irreducible summand of $N$ is of the form
$$
        W\otimes \mathbb C_{S^1},
$$
where $\mathbb C_{S^1}$ denotes a non-trivial irreducible real $S^1$-representation of
dimension $2$, and $W$ is a non-trivial irreducible real $H$-representation. We show that the case when only one-dimensional
summands $W$ occur in the decomposition of $N$ into real irreducible summands is excluded by \Cref{lemma:index_two_excl2} applied with parameters $G=S^1\times H$, $H=S^1\times H'$, $H_1=S^1$, and $H_2=H'$, where $H'$ denotes the kernel of $N$. Thus there exists $W$ with $\dim W>1$ and assumption (2) implies
$$
        \dim W\geq k.
$$
Hence
$$
        \dim (W\otimes \mathbb C_{S^1})\geq 2k.
$$
Consequently,
$$
        \dim \Sigma
        =
        \dim V
        =
        \dim V^{S^1}+\dim N
        \geq
        m+2k.
$$
This proves the assertion for $n=1$.

Now assume that the assertion has been proved for $n-1$, and consider an action of
$T^n\times H$ on a homology sphere $\Sigma$ with exactly one fixed point $x$. Write
$$
        T^n=S^1\times T^{n-1}.
$$
Let
$
        F=\Sigma^{S^1}
$. Then $F$ is a
homology sphere, and the group $T^{n-1}\times H$ acts on $F$. Moreover,
$$
        F^{T^{n-1}\times H}
        =
        \left(\Sigma^{S^1}\right)^{T^{n-1}\times H}
        =
        \Sigma^{T^n\times H}
        =
        \{x\}.
$$

We use the induced action of $T^{n-1}\times H$ on $F$. 
The inductive hypothesis applies to the effective induced action on $F$, and we
obtain
$$
        \dim F\geq m+2(n-1)k.
$$
Equivalently,
$$
        \dim V^{S^1}\geq m+2(n-1)k.
$$

It remains to estimate the normal part to $F$ at $x$. Since the original $T^n\times H$-action
is effective, the chosen $S^1$-factor acts non-trivially on the normal representation
$$
        V/V^{S^1}.
$$
As in the base case, every irreducible summand of this normal representation which is not
already excluded is of the form
$$
        W\otimes \mathbb C_{S^1},
$$
where $W$ is a non-trivial irreducible real $H$-representation. Again, by \Cref{lemma:index_two_excl2}
and assumption (2), we have $\dim W\geq k$. Therefore
$$
        \dim(V/V^{S^1})\geq 2k.
$$
Hence
$$
        \dim \Sigma
        =
        \dim V
        =
        \dim V^{S^1}+\dim(V/V^{S^1})
        \geq
        \bigl(m+2(n-1)k\bigr)+2k
        =
        m+2nk.
$$

Thus $T^n\times H$ cannot admit a one-fixed-point action on a homology sphere of dimension
less than $m+2nk$.
\end{proof}
Recall that $A_5$ admits a one fixed point action on $S^n$ whenever $n\geq 6$ and a one fixed point action on a $\mathbb{Z}$-homology sphere of dimension $3$. Moreover, the smallest possible dimension of such an action for $S_5$ on a $\mathbb{Z}$-homology sphere is $6$. We refer the reader to the following papers: \cite{BakMorimoto_equiv,BakMorimoto,Schultz1990,DeMichelis1989,Furuta1989,Morimoto1987,Morimoto1987A5,Morimoto19991A5,morimotononsol2}.
\begin{corollary}
    Let $H$ be equal to $A_5$ (resp. $S_5$). Then $H\times T^m$ does not admit one fixed point action on $n$-dimensional $\mathbb{Z}$-homology sphere, where $n<3+6m$ (resp. $n<6+8m$).
\end{corollary}
\bibliographystyle{amsplain}
\bibliography{references.bib}

\vspace{1cm} 

\bigskip
\noindent
\textit{Faculty of Mathematics and Computer Science \\
Adam Mickiewicz University \\
Uniwersytetu Poznańskiego 4 Street \\
61-614 Poznań, Poland \\
Email address:} piotr.mizerka@amu.edu.pl

\vspace{0.8cm}

\noindent
\textit{Faculty of Mathematics and Computer Science \\
Adam Mickiewicz University \\
Uniwersytetu Poznańskiego 4 Street \\
61-614 Poznań, Poland \\
Email address:} kamil.smietaniak@amu.edu.pl

\end{document}